\documentclass[preprint,sort&compress,1p]{elsarticle}
\usepackage{amssymb}
\usepackage{amsmath}
\usepackage{epstopdf}
\usepackage{topcapt}
\usepackage{color}
\usepackage{lscape}
\usepackage{psfrag}
\usepackage{amsthm}
\usepackage{url}
\usepackage{pifont}
\usepackage{natbib}
\usepackage{geometry}
\usepackage{fleqn}
\usepackage{txfonts}
\usepackage{hyperref}
\usepackage{wasysym}
\usepackage{bm}

\journal{International Journal for Numerical Methods in Fluids}

\begin{document}
\begin{frontmatter}

\title{Third-order Finite Volume\slash Finite Element Solution of the Fully Nonlinear Weakly Dispersive Serre Equations\tnoteref{label1}}
\tnotetext[label1]{The work undertaken by the first author was supported financially by an Australian National University Postgraduate Research Award}

\author[msi]{C.~Zoppou\corref{cor1}}
\ead{Christopher.Zoppou@anu.edu.au}
\cortext[cor1]{Corresponding author}

\author[msi]{J.~Pitt}
\ead{Jordan.Pitt@anu.edu.au}

\author[msi]{S.~G.~Roberts}
\ead{Stephen.Roberts@anu.edu.au}

\address[msi]{Department of Mathematics, Mathematical Sciences Institute, Australian National
University, Canberra, ACT 0200, Australia}

\begin{abstract}

The nonlinear weakly dispersive Serre equations contain higher-order dispersive terms. This includes a mixed derivative flux term which is difficult to handle numerically. The mix spatial and temporal derivative dispersive term is replaced by a combination of temporal and spatial terms. The Serre equations are re-written so that the system of equations contain homogeneous derivative terms only. The reformulated Serre equations involve the water depth and a new quantity as the conserved variables which are evolved using the finite volume method. The remaining primitive variable, the velocity is obtained by solving a second-order elliptic equation using the finite element method. To avoid the introduction of numerical dispersion that may dominate the physical dispersion, the hybrid scheme has third-order accuracy. Using analytical solutions, laboratory flume data and by simulating the dam-break problem, the proposed scheme is shown to be accurate, simple to implement and stable for a range of problems, including discontinuous flows.

\end{abstract}

\begin{keyword}
dispersive waves \sep conservation laws \sep Serre equation \sep finite volume \sep finite elements

\MSC[2010] 76B15 \sep 35L65 \sep 65M08

\end{keyword}
\end{frontmatter}


\section{Introduction} \label{intro}

Rapidly-varying free surface flows are characterized by large surface gradients. These gradients produce vertical accelerations of fluid particles resulting in a non-hydrostatic pressure distribution. Equations that assume that the flow has a non-hydrostatic pressure distribution contain third-order dispersive terms. A system of equations that contain dispersive terms are the Serre equations. The Serre equations are fully nonlinear weakly dispersive equations. They are applicable to waves where the water depth $h_0 \ll L$ is much smaller than the horizontal wave length, $L$ and up to wave breaking\cite{Barthelemy-E-2004-315}. Bonneton \emph{et al.}\cite{Bonneton-etal-2011-1479,Bonneton-etal-2011-589} consider the Serre equations as the most appropriate system for modelling highly nonlinear weakly dispersive waves at the shoreline.

A major difficulty with solving equations that contain dispersive terms, is that the dispersive terms usually contain a mix derivative term\cite{Dias-Milewski-2010}. By replacing the mix derivative term in the flux term by a combination of  temporal and spatial derivative terms, the Serre equations can be written in conservation law form, where the system of homogeneous equations contains a new conserved quantity and its corresponding flux term. The conserved quantities are evolved using a standard scheme for solving conservative laws. The remaining primitive variable is obtained by solving a second-order elliptic equation. A finite volume\slash finite element technique is proposed for the solution of the fully nonlinear and weakly dispersive Serre equations without the need for iteration or operator splitting necessary for dealing with the mixed derivative term.

It is well known that odd-order schemes introduce numerical dispersion and even-order schemes numerical dispersion.  First-order and second-order schemes were developed by Zoppou and Roberts\cite{Zoppou-Roberts-2013}, using the methodology described in this paper, to solve the Serre equation. Both schemes produced dispersive waves that accurately predict the arrival of the initial wave and its amplitude. In addition, the phase of the predicted dispersive waves is very close to the recorded wave profile. However, the diffusion introduced by the first-order scheme rapidly dampens trailing dispersive waves. The second-order model slightly overestimates the amplitude of the dispersive waves. This could be a result of the use of a second-order scheme which introduced numerical dispersion. To avoid the introduction of numerical dispersion and excessive diffusion, a third-order finite volume\slash finite element scheme was developed.

The performance of the proposed third-order finite volume/finite element scheme for solving the conservative form of the Serre equations is evaluated using an analytical solution to the Serre equations, laboratory flume data and by simulating the dam-break problem. With the exception of the analytical solution, which is smooth, the remaining problems involve the simulation of flows with steep gradients that produce dispersive waves.

This paper is organized with the derivation of the standard Serre equation in Section 2 followed by the derivation of the alternative form of the Serre equation in terms of the new conservative variable. The properties of the linearized form of the Serre equation are also examined in Section 2. The third-order implementation of the proposed scheme is described in detail in Section 3 and the fully discretized scheme summarized in Section 4. In Section 5, the numerical scheme is validated using an analytical solution and laboratory flume data. The stability of the proposed scheme is demonstrated by simulating the dam-break problem. Finally, the performance of the numerical scheme is discussed in Section 6.

\section{Serre Equations}

\begin{figure}[htb]
\begin{center}
\includegraphics[width=7.0cm]{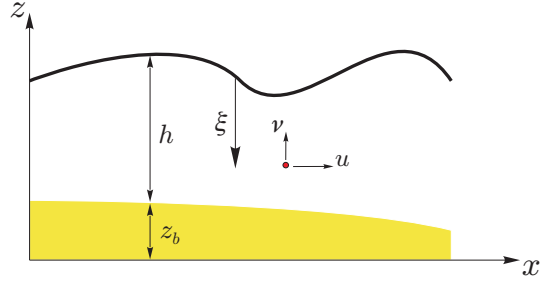}
\end{center}
\caption{The notation used for one-dimensional flow governed by the Serre equation.}
\label{fig:Notation}
\end{figure}

For an invicid incompressible fluid with constant density, $\rho$ the conservation of mass and momentum are given by the Euler equations
\begin{subequations}
\begin{equation}
\nabla \cdot \textbf{u} =  0,
\label{eq:Euler_continuity}
\end{equation}
and
\begin{equation}
\rho \frac{D\textbf{u}}{Dt} = - \nabla p + \rho \textbf{g}
\label{eq:Euler_momentum}
\end{equation}
\label{eq:Euler_governing_equations}
\end{subequations}
where in two planar dimensions, $\textbf{x} = (x,z)$, a fluid particle at depth  $\xi = z - h - z_b$ below the water surface, where $h(x,t)$ is the water depth and $z_b(x)$ the bed elevation, see Figure \ref{fig:Notation}, is subject to fluid pressure, $p(\textbf{x},t)$ and  gravitational acceleration, $\textbf{g} = (0,g)^T$, has a velocity $\textbf{u} = (u(\textbf{x},t),v(\textbf{x},t))$,  where $u(\textbf{x},t)$ is the velocity in the $x$-coordinate and $v(\textbf{x},t)$ is the velocity in the $z$-coordinate and $t$ is time.
In addition to the above equations, a number of boundary conditions must be satisfied. These are;
\begin{subequations}\label{eq:Seabra_all}
\begin{enumerate}[(a)]
\item the kinematic condition at the free surface $(z = h + z_b)$,
\begin{equation*}
v|_{h+z_b} = \dfrac{\partial h}{\partial t} + u \dfrac{\partial (h + z_b)}{\partial x},
\end{equation*}
\item the kinematic condition at the bed $(z = z_b)$,
\begin{equation*}
v|_{z_b} = u \dfrac{\partial z_b}{\partial x}
\end{equation*}
\item the dynamic condition at the surface $(z = h + z_b)$
\begin{equation*}
p(\xi = 0) = p_a.
\end{equation*}
\end{enumerate}
\end{subequations}
which is the atmospheric pressure at the water surface, usually taken to be $p_a = 0$.

The Serre equations assume that the point velocity in the $x$-direction is uniform over the water depth, so that $u(x,z,t) = \bar{u}(x,t)$ with
\begin{equation*}
\bar{u}(x,t) = \dfrac{1}{h} \int_{z_b}^{h+z_b} u(x,z,t) \,dz.
\end{equation*}

From \eqref{eq:Euler_continuity} it follows that the vertical velocity at any depth $z - z_b$ is given by
\begin{equation*}
v |_z = -(z - z_b) \dfrac{\partial \bar{u}}{\partial x}
\end{equation*}
for a horizontal bed. The vertical velocity is a linear function of the water depth.

Integrating the point quantities in \eqref{eq:Euler_momentum} over the flow depth $z_b$ to $h+z_b$, and satisfying \eqref{eq:Seabra_all} produces the one-dimensional equations
\begin{subequations}\label{eq:Serre_nonconservative_form}
\begin{equation}
\dfrac{\partial h}{\partial t} + \bar{u} \dfrac{\partial
h}{\partial x} +  h \dfrac{\partial
\bar{u}}{\partial x} = 0 %
\label{eq:Boussinsq_continuity}
\end{equation}
and
\begin{equation}
\dfrac{\partial \bar{u}}{\partial t} + \bar{u} \dfrac{\partial \bar{u}}{\partial x} + \dfrac{1}{h}\dfrac{\partial}{\partial x} \left [  \dfrac{gh^2}{2} + \dfrac{h^3}{3} \left ( \dfrac{\partial \bar{u}}{\partial x} \dfrac{\partial \bar{u}}{\partial x} - \bar{u} \dfrac{\partial^2\bar{u}}{\partial x^2} - \dfrac{\partial^2\bar{u}}{\partial x \partial t} \right ) \right ] = 0
\label{eq:Boussinsq_momentum}
\end{equation}
\end{subequations}
where $h$ and $\bar{u}$ are the primitive variables.

The pressure distribution in the water column is given by
\begin{equation}
p|_\xi = p_a + \rho g \xi + \dfrac{\rho}{2} \xi ( 2h - \xi ) \left ( \dfrac{\partial \bar{u}}{\partial x} \dfrac{\partial \bar{u}}{\partial x} - \bar{u} \dfrac{\partial^2 \bar{u}}{\partial x^2} -  \dfrac{\partial^2 \bar{u}}{\partial x \partial t} \right ).
\end{equation}
Multiplying \eqref{eq:Boussinsq_momentum} by $h$, adding \eqref{eq:Boussinsq_continuity} pre-multiplied by $\bar{u}$ and making use of \eqref{eq:Boussinsq_continuity} to obtain;
\begin{subequations}\label{eq:Serre_conservative_form}
\begin{equation}
\dfrac{\partial h}{\partial t} + \dfrac{\partial (\bar{u}h)}{\partial x} = 0
\label{eq:Serre_continuity}
\end{equation}
and
\begin{equation}
\dfrac{\partial (\bar{u}h)}{\partial t} + \dfrac{\partial}{\partial x} \left ( \bar{u}^2h + \dfrac{gh^2}{2}\right ) + \dfrac{\partial}{\partial x} \left [  \dfrac{h^3}{3} \left ( \dfrac{\partial \bar{u}}{\partial x} \dfrac{\partial \bar{u}}{\partial x} - \bar{u} \dfrac{\partial^2\bar{u}}{\partial x^2} - \dfrac{\partial^2 \bar{u}}{\partial x \partial t} \right ) \right ]  = 0
\label{eq:Serre_momentum}
\end{equation}
\end{subequations}
which is written in terms of the conservative variables, $h$ and $\bar{u}h$. The continuity equation is exact because it is based on the depth-averaged velocity, which makes no assumpltion on the distribution of $u(x,z,t)$ with depth. However, the momentum equation relies on the assumption that $u(x,z,t)$ is uniform with depth.

The terms in the square parenthesis are the dispersive terms which contain high order spatial derivative terms and a mixed spatial and temporal derivative term. Ignoring all the dispersive terms in \eqref{eq:Serre_momentum} results in the well known nonlinear shallow water wave equations, where the pressure distribution is hydrostatic,  $p(\xi) = \rho g \xi$.

Equation \eqref{eq:Serre_conservative_form} are known as the Serre equations\cite{Serre-F-1953-857,Seabra-Santos-etal-1987-117,Carter-Cienfuegos-2010-259}, they retain full nonlinearity in the dispersive terms\cite{El-etal-2006}. They have been derived by Serre\cite{Serre-F-1953-857}, Su and Gardner\cite{Su-Gardener-1969-536} and Seabra-Santos \emph{et al.}\cite{Seabra-Santos-etal-1987-117} and are equivalent to the depth averaged Green and Naghadi\cite{Green-Naghdi-1976-237} equations. They are considered to be good approximations to the full Euler equations up to a wave breaking\cite{Bonneton-etal-2011-1479,Bonneton-etal-2011-589}.

\subsection{Alternative Conservative Form of the Sere Equations}

The flux term in the momentum equation, \eqref{eq:Serre_momentum} contains a mixed derivative term which is difficult to treat numerically. It is possible to replace the mix spatial and temporal derivative term  by a combination of spatial and temporal derivative terms.

Consider
\begin{equation*}
\dfrac{\partial^2}{\partial x \partial t} \left ( \dfrac{h^3}{3} \dfrac{\partial \bar{u}}{\partial x} \right ) = \dfrac{\partial }{\partial t} \left ( h^2 \dfrac{\partial h}{\partial x} \dfrac{\partial \bar{u}}{\partial x} + \dfrac{h^3}{3} \dfrac{\partial^2 \bar{u}}{\partial x^2} \right ) =
\dfrac{\partial }{\partial x} \left ( h^2 \dfrac{\partial h}{\partial t} \dfrac{\partial \bar{u}}{\partial x} + \dfrac{h^3}{3} \dfrac{\partial^2 \bar{u}}{\partial x \partial t} \right ).
\end{equation*}
Rearranging then
\begin{equation*}
\dfrac{\partial }{\partial x} \left ( \dfrac{h^3}{3} \dfrac{\partial^2 \bar{u}}{\partial x \partial t} \right ) = \dfrac{\partial }{\partial t} \left ( h^2 \dfrac{\partial h}{\partial x} \dfrac{\partial \bar{u}}{\partial x} + \dfrac{h^3}{3} \dfrac{\partial^2 \bar{u}}{\partial x^2} \right ) - \dfrac{\partial }{\partial x} \left ( h^2 \dfrac{\partial h}{\partial t} \dfrac{\partial \bar{u}}{\partial x}\right ).
\end{equation*}
Making use of the continuity equation, \eqref{eq:Serre_continuity}
\begin{equation*}
\dfrac{\partial }{\partial x} \left ( \dfrac{h^3}{3} \dfrac{\partial^2 \bar{u}}{\partial x \partial t} \right ) = \dfrac{\partial }{\partial t} \left ( h^2 \dfrac{\partial h}{\partial x} \dfrac{\partial \bar{u}}{\partial x} + \dfrac{h^3}{3} \dfrac{\partial^2 \bar{u}}{\partial x^2} \right ) + \dfrac{\partial }{\partial x} \left [ h^2 \dfrac{\partial \bar{u}}{\partial x} \left ( \bar{u}\dfrac{\partial h}{\partial x} + h \dfrac{\partial \bar{u}}{\partial x}  \right ) \right ]
\end{equation*}
and the momentum equation, \eqref{eq:Serre_momentum} becomes
\begin{equation*}
\dfrac{\partial }{\partial t} \left ( \bar{u}h -  h^2 \dfrac{\partial h}{\partial x} \dfrac{\partial \bar{u}}{\partial x} - \dfrac{h^3}{3} \dfrac{\partial^2 \bar{u}}{\partial x^2}  \right ) + \dfrac{\partial}{\partial x} \left ( \bar{u}^2h + \dfrac{gh^2}{2} - \bar{u}h^2 \dfrac{\partial h}{\partial x}\dfrac{\partial \bar{u}}{\partial x} - \dfrac{\bar{u}h^3}{3} \dfrac{\partial^2 \bar{u}}{\partial x^2} - \dfrac{2h^3}{3}\dfrac{\partial \bar{u}}{\partial x}\dfrac{\partial \bar{u}}{\partial x} \right ) \nonumber  = 0.
\end{equation*}
The momentum equation can be written in terms of a new conservative form as
\begin{equation*}
\dfrac{\partial G}{\partial t} + \dfrac{\partial }{\partial x} \left ( G\bar{u} + \dfrac{gh^2}{2} - \dfrac{2 h^3}{3} \dfrac{\partial \bar{u}}{\partial x} \dfrac{\partial \bar{u}}{\partial x} \right ) = 0
\end{equation*}
where the new conserved quantity, $G$ is given by the second-order elliptic equation
\begin{equation}
G = \bar{u}h - \dfrac{\partial }{\partial x} \left ( \dfrac{h^3}{3} \dfrac{\partial \bar{u}}{\partial x} \right ) \label{eq:G_FE}.
\end{equation}
The temporal derivative in the momentum equation has been eliminated from the flux term. In contrast to \eqref{eq:Serre_conservative_form}, the flux term now contains spatial derivatives only.

For the remaining primitive variable, $\bar{u}$, if the data is square integrable in a rectangular domain, then from the regularity theorem of elliptic partial differential equations\cite{Evans-L-1997} $\bar{u} \in H^2$ in Sobolev space of square integrable second-derivatives. The primitive variable, $\bar{u}$ will be smooth.

The alternative form of the Serre equations can be written in vector form as
\begin{subequations}
\label{eq:Serre_Vector}
\begin{equation}
\dfrac{\partial \bm{q}(x,t)}{\partial t} + \dfrac{\partial \bm{F}(\bm{q}(x,t))}{\partial x} = 0.
\end{equation}
where the vector of state variables
\begin{equation}\label{eq:Serre_conserved_quantity}
\bm{q}(x,t) = \left[ \begin{array}{c} h \\ G
\end{array} \right],
\end{equation}
and
\begin{equation}\label{eq:Serre_flux}
\bm{F}(\bm{q}(x,t)) = \left [ \begin{array}{c} f(1) \\
f(2) \end{array} \right ] = \left[ \begin{array}{c} \bar{u}h \\
G\bar{u} + \dfrac{gh^2}{2} - \dfrac{2h^3}{3} \dfrac{\partial
\bar{u}}{\partial x} \dfrac{\partial \bar{u}}{\partial x}
\end{array} \right].
\end{equation}
\end{subequations}

\subsubsection{Properties of the Linearized Serre equations}

Although they are evolution-type equations, the Serre equations are neither hyperbolic or parabolic and do not have any Riemann invariants. However, it is possible to establish some properties of the Serre equations by examining the behaviour of harmonic waves of the form
\begin{gather}
\label{eq:Fourier_components}
h(x,t) = A e^{i(kx - \omega t)} \quad \text{and} \quad u(x,t) = U e^{i(kx - \omega t)}
\end{gather}
where $A$ and $U$ are unknown coefficients, $\omega$ is the frequency, $k = 1/L$ is the  wave number and $i = \sqrt{-1}$
in the linearized Serre equations
\begin{subequations}
\label{eq:linearized_Serre}
\begin{gather}
\label{eq:linearized_Serre_continuity}
\dfrac{\partial h_1}{\partial  t} + h_0 \dfrac{\partial u_1}{\partial x} +  u_0 \dfrac{\partial h_1}{\partial x} = 0
\end{gather}
and
\begin{gather}
\label{eq:linearized_Serre_momentum}
\dfrac{\partial u_1}{\partial t} +  g \dfrac{\partial h_1}{\partial x} + u_0 \dfrac{\partial u_1}{\partial x} -  \dfrac{h_0^2}{3} \left ( u_0 \dfrac{\partial^3 u_1}{\partial x^3}  + \dfrac{\partial^3 u_1}{\partial x^2 \partial t} \right )  = 0 .
\end{gather}
\end{subequations}

Substituting \eqref{eq:Fourier_components} into \eqref{eq:linearized_Serre}, the linearized equations become
\begin{subequations}
\begin{gather*}
-A \omega + u_0 A k + h_0 U k = 0
\end{gather*}
and
\begin{gather*}
- U \omega +  g A k +  u_0 U k - \dfrac{1}{3} h_0^2 U  \omega k^2 + \dfrac{1}{3} h_0^2 u_0 U  k^3 = 0.
\end{gather*}
\end{subequations}
For a non-trivial solution
\begin{gather*}
\left|
  \begin{array}{cc}
     -\omega + u_0 k & h_0 k \\
    g k & -\omega + u_0 k - \dfrac{1}{3} h_0^2 \omega k^2 + \dfrac{1}{3} h_0^2 u_0 k^3 \\
  \end{array}
\right| = 0
\end{gather*}
or
\begin{gather*}
\omega_{1,2} = u_0 k \pm k \sqrt{g h_0} \sqrt{\dfrac{3}{\mu^2 + 3}} \nonumber
\end{gather*}
where $\mu = h_0k$ is the frequency dispersion.

In this case the dispersive terms have no effect on $u_0$, only on the celerity of a small disturbance.

This compares with the frequency, $\omega_{1,2} = (\bar{u}_0 \pm \sqrt{g h_0}) k$, for the shallow water wave equations. As $\mu \rightarrow 0$, the frequency for the Serre equations are identical to that of the shallow water wave equations. When $\mu \rightarrow \infty$, $\omega_{1,2} = \bar{u}_0$. Therefore, the frequency for the Serre equation are bounded by the wave frequency of the shallow water wave equations.

For non-dispersive waves, the phase velocity, $\upsilon_p = \text{Re}(\omega)/k$ is identical to the group velocity $\upsilon_g = d\text{Re}(\omega)/dk$.  This is not the case for the Serre equations, where the phase speed is
\begin{subequations}
\begin{gather*}
\upsilon_p = u_0 \pm \sqrt{gh_0}\sqrt{\dfrac{3}{\mu^2 + 3}}
\end{gather*}
and the group velocity is
\begin{gather*}
\upsilon_g = u_0 \pm \sqrt{gh_0}\left ( \sqrt{\dfrac{3}{\mu^2 + 3}} \mp \mu^2 \sqrt{\dfrac{3}{(\mu^2 + 3)^3}} \right ) \neq \upsilon_p.
\end{gather*}
\end{subequations}
Both are dependent on the wave number. Since the group speed is slower than the phase speed then the Serre equations describe dispersive waves.

\section{Numerical Scheme}\label{numerical_scheme}

In a finite volume scheme the cell average values at the nodes, $\bar{q}_j$ are updated by integrating \eqref{eq:Serre_Vector} over the $j$th cell $I_j = [x_{j-1/2},x_{j+1/2}]$ to obtain a semi-discrete scheme
\begin{equation}
\label{eq:semi-discrete}
\dfrac{d\bar{q}_j(t)}{dt} = \mathcal{L}(\bar{q}(x,t)) = -\dfrac{f(q(x_{j+1/2},t)) - f(q(x_{j-1/2},t))}{\Delta x}.
\end{equation}
where $\Delta x = x_{j+1/2} - x_{j-1/2}$ is assumed to be constant, $x_j = (x_{j+1/2} + x_{j-1/2})/2$ and
\begin{equation*}
\bar{q}_j(t) = \dfrac{1}{\Delta x} \int_{x_{j-1/2}}^{x_{j+1/2}} q(x,t) \, dx
\end{equation*}
is the average value of the state variable $q(x,t)$ in $I_j$ at time $t$, which ensures that mass is conserved in each cell and the integral formulation admits shocks in the solution.

In the discrete forms of \eqref{eq:Serre_Vector}, $f(q(x_{j \pm 1/2},t)) = f_{j \pm 1/2}(\bar{q}_{j-1}(t), \dots, \bar{q}_{j+1}(t)) = F_{j \pm 1/2}$ represents the numerical approximation of the physical flux $f(q(x,t))$ across the boundary of cell $j$, at $x_{j \pm 1/2}$ at time $t$.

The flux , $F_{j+1/2}$ is a function of the left and right extrapolated state values $q^+_{j+1/2}$ and $q^-_{j+1/2}$, obtained from piecewise polynomials, $P_j(x_{j+1/2})$ and $P_{j+1}(x_{j+1/2})$, respectively passing through consecutive values of $\bar{q}_j$. Therefore,
\begin{equation*}
F_{j+1/2} = f_{j+1/2}(q_{j+1/2}^+,q_{j+1/2}^-).
\end{equation*}

The reconstruction will usually produce two different values for $q(x_{j+1/2})$. Generally, there will be a discontinuity in the state variables at $x_{j\pm 1/2}$.

The flux of material across the interface of a cell is estimated by solving the Riemann problem, defined by the initial value problem
\begin{equation*}
q(x_{j+1/2}) = \left \lbrace \begin{matrix} q_{j+1/2}^+ \;\; \text{if} \;\; x < x_{j+1/2} \\
q_{j+1/2}^- \;\; \text{if} \;\; x > x_{j+1/2}.
\end{matrix}
\right .
\end{equation*}
Once the intercell flux has been established, the cell average values can be updated by evolving the solution over a single time step by solving the semi-discrete system, \eqref{eq:semi-discrete} using an ordinary differential equation solver.  The overall accuracy of the numerical scheme is dependent on the accuracy of the reconstruction method and the order of accuracy of the scheme used to evolve the solution in time.

\subsection{Inter-cell Flux Evaluation}

The numerical approximation of the physical flux $f(q(x,t))$ across the boundary of a cell, $F_{j+1/2}$ is given by the explicit upwind central scheme proposed by Kurganov \emph{et al.}\cite{Kurganov-etal-2001-707} as
\begin{equation}
\label{eq:HLL_flux}
F_{j+1/2} = \dfrac{a_{j+1/2}^+ f(q^-_{j+1/2}) -a_{j+1/2}^- f(q^+_{j+1/2})}{a^+_{j+1/2} - a^-_{j+1/2}}  + \dfrac{a_{j+1/2}^+ \, a_{j+1/2}^-}{a_{j+1/2}^+ - a_{j+1/2}^-} \left [ q^+_{j+1/2} - q^-_{j+1/2} \right ].
\end{equation}

At the interface of a cell, $x_{j \pm 1/2}$ a discontinuity in the state variable will propagate with right- and left-sided local speeds, which are estimated by
\begin{subequations}
\begin{equation*}
a^+_{j+1/2} = \text{max}\left [ \lambda_2(q_{j+1/2}^-), \, \lambda_2(q_{j+1/2}^+), \, 0 \right ],
\end{equation*}
and
\begin{equation*}
a^-_{j+1/2} = \text{min}\left [ \lambda_1 (q_{j+1/2}^-), \, \lambda_1 (q_{j+1/2}^+), \,0 \right ]
\end{equation*}
\end{subequations}
where $\lambda_1$ and $\lambda_2$ are the smallest and largest eigenvalues, respectively  of the Jacobian $\partial f(\bar{q})/\partial \bar{q}$ which correspond to the phase speeds.

It has been established that the eigenvalues are bounded for the Serre equations and that the eigenvalues,  $\lambda_1 \simeq \bar{u} + \sqrt{gh}$, $\lambda_2 \simeq \bar{u} - \sqrt{gh}$ for the shallow water wave equation can be used as estimates of the upper and lower bounds of the eigenvalues for the Serre equations.

There is a restriction on the computational time-step that can be used in all explicit schemes. Stability is satisfied when the time step $\Delta t$ satisfies the \emph{Courant-Friedrichs-Lewy}, ($CFL$) criteria\cite{Harten-etal-1983-357}
\begin{equation*}
\Delta t < \dfrac{\Delta x}{2 \text{max}(|\lambda_i|)}  \quad \forall
\,\, i.
\end{equation*}

\subsection{Reconstruction}

To achieve third-order, $\mathcal{O}(\Delta x^3)$ accuracy, it is sufficient to consider quadratic piecewise polynomial reconstruction in each cell using
\begin{equation*}
P_j(x) = a_j(x - x_j)^2 + b_j(x - x_j) + c_j \quad x \in [x_{j-1/2},x_{j+1/2}].
\end{equation*}
centered at the computational node $j$. Over each cell, the cell average
\begin{equation*}
q^n_j = \dfrac{1}{\Delta x} \int_{x_{j-1/2}}^{x_{j+1/2}} P_j(x) \, dx.
\end{equation*}
The mass conserving quadratic is given by
\begin{equation}
P_j(x) = \bar{q}_j + \dfrac{(x - x_j)}{\Delta x} \dfrac{\bar{q}_{j+1} - \bar{q}_{j-1}}{2}
 + 3 \kappa \left [ (x - x_j)^2 - \dfrac{\Delta x^2}{12} \right ] \dfrac{\bar{q}_{j+1} - 2 \bar{q}_j + \bar{q}_{j-1}}{2 \Delta x^2}
\label{eq:Hirsh_21.1.1}
\end{equation}
where $-1 \le \kappa \le 1$. Only when $\kappa = 1/3$ is the reconstruction third-order.  Any other value of $\kappa$ results in linear reconstruction and a second-order scheme.

The reconstructed cell edge values at $x_{j\pm1/2}$ are estimated from \eqref{eq:Hirsh_21.1.1}. Consider the cell interface at $x_{j+1/2}$, then \eqref{eq:Hirsh_21.1.1} becomes
\begin{equation}
\label{eq:Hirsh_21.1.5a}
q^-_{j +1/2} = \bar{q}_j + \dfrac{1}{4} ( 1 - \kappa )(\bar{q}_j - \bar{q}_{j-1}) + \dfrac{1}{4} ( 1 + \kappa )(\bar{q}_{j+1} - \bar{q}_j )
\end{equation}
which makes use of $\kappa/2 = -1/4 + \kappa/4 + 1/4 + \kappa/4$.

Similarly at the left boundary of the cell, at $x_{j-1/2}$
\begin{equation*}
q^+_{j-1/2} = \bar{q}_j - \dfrac{1}{4} ( 1 + \kappa )(\bar{q}_j - \bar{q}_{j-1}) - \dfrac{1}{4} ( 1 - \kappa )(\bar{q}_{j+1} - \bar{q}_j )
\end{equation*}

The interpolated values must be limited so that no new extremum is introduced.  Consider \eqref{eq:Hirsh_21.1.5a} and making use of $\bar{q}_{j+1} - \bar{q}_{j-1} = \bar{q}_{j+1} - \bar{q}_j + ( \bar{q}_j - \bar{q}_{j-1} )$, then
\begin{equation*}\label{eq:second_order_general}
q^-_{j+1/2} =  \bar{q}_j + \dfrac{1}{2} \left [ \left ( \dfrac{1}{2} + \dfrac{\kappa}{2} \right ) (\bar{q}_{j+1} - \bar{q}_j ) + \left ( \dfrac{1}{2} - \dfrac{\kappa}{2} \right )(\bar{q}_j - \bar{q}_{j-1} ) \right ].
\end{equation*}
Dividing the term in the square brackets by $\bar{q}_{j} - \bar{q}_{j-1}$ then
\begin{equation*}
q^-_{j+1/2} = \bar{q}_j + \dfrac{1}{2} \left [ \left ( \dfrac{1}{2} + \dfrac{\kappa}{2} \right ) \dfrac{\bar{q}_{j+1} - \bar{q}_j}{\bar{q}_j - \bar{q}_{j-1}} + \left ( \dfrac{1}{2} - \dfrac{\kappa}{2} \right ) \right ] (\bar{q}_j - \bar{q}_{j-1} ). \end{equation*}
Let $r_j = (\bar{q}_{j+1} - \bar{q}_j)/(\bar{q}_j - \bar{q}_{j-1})$, then
\begin{equation*}
q^-_{j+1/2} = \bar{q}_j + \dfrac{1}{2} \left [ \left ( \dfrac{1}{2} + \dfrac{\kappa}{2} \right ) r_j + \left ( \dfrac{1}{2} - \dfrac{\kappa}{2} \right ) \right ] (\bar{q}_j - \bar{q}_{j-1} ).
\end{equation*}
When $\kappa = 1/3$ then
\begin{eqnarray}
q^-_{j+1/2} &=&  \bar{q}_j + \dfrac{1}{2} \left ( \dfrac{2}{3} r_j + \dfrac{1}{3} \right ) ( \bar{q}_j - \bar{q}_{j-1} ) \nonumber \\
&=& \bar{q}_j + \dfrac{1}{2} \phi^-(r_j) (\bar{q}_j - \bar{q}_{j-1} ) \label{eq:limiter_u_minus}
\end{eqnarray}
where
\begin{equation*}\label{eq:limiter_phi_minus}
\phi^-(r_j) = \dfrac{2}{3} r_j + \dfrac{1}{3}
\end{equation*}
is a nonlinear limiter which can be used to prevent unwanted oscillations and ensures that the results are physical (bounded) and therefore stable.

Equation \eqref{eq:limiter_u_minus} forms the basis of the third-order symmetrical Koren limiter\cite{Koren-B-1993} given by
\begin{equation*}
\phi^-(r_j) = \text{max}\left [0,\text{min}(2r_j, (1 + 2r_j)/3,2) \right ] \label{eq:koren_limiter}
\end{equation*}
where $\lim_{r_j \rightarrow \infty} \phi^+(r_j) = 2$ is used in the proposed model. This limiter ensures that the scheme remains \emph{TVD}\cite{Sweby-P-1984-995} and third-order away from extrema, where $r_j \le 0$ and $\phi(r_j) = 0$.  This occurs when the gradient changes sign indicating that an extrema has been encountered within the cell.  The reconstruction reverts to a piecewise constant reconstruction. In smooth regions, $r_j \rightarrow 1$ and $\phi(r_j) \rightarrow 1$ and the reconstruction is third-order.

Similarly for the reconstructed values at the left boundary of the cell, the limited values are given by
\begin{equation}
\label{eq:limiter_u_plus}
q^+_{j+1/2} = \bar{q}_j - \dfrac{1}{2} \phi^+(r_j) (\bar{q}_j - \bar{q}_{j-1} ) \nonumber
\end{equation}
where
\begin{equation*}
\phi^+(r_j) = \dfrac{2}{3} + \dfrac{1}{3}r_j
\end{equation*}
and the corresponding nonlinear limiter is
\begin{equation*}
\phi^+(r_j) = \text{max}\left [0,\text{min}(2r_j, (2 + r_j)/3,2) \right ].
\end{equation*}

\subsection{Nodal Values}

The point values for the conservative quantities, $q_j$ are estimated from the cell averages, $\bar{q}(x_j)$ using quadratic interpolation, \eqref{eq:Hirsh_21.1.1} so that
\begin{align}
q_j =  \dfrac{- \bar{q}(x_{j+1}) + 26 \bar{q}(x_j) - \bar{q}(x_{j-1})}{24}.
\label{eq:cell_to_node}
\end{align}
Written for all computational nodes results in a tridiagonal matrix, $\mathcal{M}$ which can be solved for the cell averages, $\bar{q}(x_j)$ given the nodal values, $q_j$.

\subsection{Time Integration}

Time integration of the semi-discrete system \eqref{eq:semi-discrete} is performed using Strong Stability Preserving (\emph{SSP}) Runge-Kutta schemes. SSP schemes involve a convex combination of first-order forward Euler steps that preserve the \emph{TVD} properties of the Euler scheme\cite{Shu-Osher-1988-439,Gottlieb-etal-2009-251}.

A third-order three-stage \emph{SSP} Runge-Kutta scheme is given by\cite{Shu-Osher-1988-439,MacDonald-etal-2008-89}
\begin{subequations}
\begin{equation}
\bar{q}^{(1)}_j = \bar{q}^n_j + \Delta t \mathcal{L}(t_n,\bar{q}^n_j),
\end{equation}
\begin{equation}
\bar{q}^{(2)}_j = \bar{q}^{(1)}_j + \Delta t \mathcal{L} \left ( t_n+\Delta t,\bar{q}^{(1)}_j \right ),
\end{equation}
\begin{equation}
\bar{q}^{(3)}_j = \dfrac{3}{4}\bar{q}^n_j + \dfrac{1}{4}\bar{q}^{(2)}_j,
\end{equation}
\begin{equation}
\bar{q}^{(4)}_j = \bar{q}^{(3)}_j + \Delta t \mathcal{L} \left ( t_n+\dfrac{\Delta t}{2},\bar{q}^{(3)}_j \right )
\end{equation}
and
\begin{equation}
\bar{q}^{n+1}_j = \dfrac{1}{3}\bar{q}^n_j + \dfrac{2}{3}\bar{q}^{(4)}_j
\end{equation}
\label{eq:RK_TVD_shu_Osher0}
\end{subequations}
where
\begin{equation}
\label{eq:Flux_difference}
\mathcal{L}(t_n,\bar{q}_j) = - \dfrac{\Delta t}{\Delta x} \left ( F_{j+1/2} - F_{j-1/2} \right )
\end{equation}
is the discrete form of \eqref{eq:semi-discrete} and It is also subject to the time restriction, $Cr = 1$.

It is relatively straight forward to establish the boundary information required for the third-order Runge-Kutta \emph{TVD} scheme by examining equation \eqref{eq:RK_TVD_shu_Osher0}. In this case boundary information is requires at the time levels; $t_n$, $t_n + \Delta t/2$, $t_n + \Delta t$, $t_n + 3\Delta t/2$ and $t_n + 2\Delta t$.

Combined with quadratic reconstruction, the resulting numerical scheme is theoretically $O(\Delta x^3, \Delta t^3)$ accurate.

\section{Fully Discrete System}

The reconstructed point values, $\bar{u}^\pm_{j+1/2}$, $h^\pm_{j+1/2}$ and $G^\pm_{j+1/2}$ at each cell interface are calculated from the cell averaged values using \eqref{eq:limiter_u_minus} and \eqref{eq:limiter_u_plus}. These are used to estimate the left and right intercell flux using a finite-difference discretization of \eqref{eq:Serre_flux}. The approximate Riemann solver \eqref{eq:HLL_flux} is used to obtain the intercell flux required in \eqref{eq:Flux_difference}. Using an \emph{SSP} Runge-Kutta scheme, the solution is advanced in time providing updated values for the cell average conservative variables, $\bar{h}^{n+1}_j$ and $\bar{G}_j^{n+1}$. The nodal values for the remaining primitive variable $\bar{u}_j^{n+1}$ is obtained by solving the second-order elliptic equation, \eqref{eq:G_FE} given point values for $h^{n+1}_j$ and $G_j^{n+1}$, obtained using \eqref{eq:cell_to_node} and the cell average values, using finite elements..

The third-order finite difference discretizations of $f(2)$ in \eqref{eq:Serre_flux} and the third-order finite element solution of \eqref{eq:G_FE} are described in detail in the following Section.

\subsection{Discretization of the Instantaneous Flux, $f$}

The spatial derivatives in the local flux term, \eqref{eq:Serre_flux} are evaluated using a quadratic polynomial fitted through the values, $u_j$ and $u_{j \pm 1/2}$, which were obtained from the finite element solution of the second-order elliptic equation, \eqref{eq:G_FE} . The dircretized  second component of the local flux can be written as
\begin{subequations}
\label{eq:local_flux_discrete}
\begin{equation}
f(2)(q^+_{j-1/2}) = (G\bar{u})^+_{j-1/2} + \dfrac{g{h^+_{j-1/2}}^2}{2} - \dfrac{2{h^+_{j-1/2}}^3}{3\Delta x^2} \left ( -u^+_{j+1/2} + 4u^+_{j} - 3u^+_{j-1/2} \right )^2,
\end{equation}
\begin{equation}
f(2)(q^-_{j+1/2}) = (G\bar{u})^-_{j+1/2} + \dfrac{g{h^-_{j+1/2}}^2}{2} - \dfrac{2{h^-_{j+1/2}}^3}{3\Delta x^2} \left ( 3u^-_{j+1/2} - 4u^-_j + u^-_{j-1/2} \right )^2
\end{equation}
and for the first component
\begin{equation}
f(1)(u^\pm_{j+1/2}) = (\bar{u}h)^\pm_{j+1/2}.
\end{equation}
\end{subequations}

\subsection{Finite Element Solution of the Elliptic Equation}

In finite elements, the second-order elliptic equation \eqref{eq:G_FE} becomes
\begin{equation*}
\int_a^b G \upsilon \,\, dx = \int_a^b \bar{u}h \upsilon \,\, dx - \int_a^b \dfrac{\partial}{\partial x} \left ( \dfrac{h^3}{3} \dfrac{\partial \bar{u}}{\partial x} \right ) \upsilon \,\, dx
\end{equation*}
where $x \in [a,b]$ is the computational domain and $\upsilon \in H_0^1(a,b)$ is a test function.

Integrating the last term by parts, then the weak form of the weighted residual equation is given by
\begin{equation*}
\int_a^b G \upsilon \,\, dx =  \int_a^b \bar{u}h \upsilon \,\, dx - \left . \upsilon \dfrac{h^3}{3} \dfrac{\partial \bar{u}}{\partial x} \right |_a^b + \int_a^b \dfrac{h^3}{3} \dfrac{\partial \bar{u}}{\partial x} \dfrac{d
\upsilon}{dx} \,\, dx
\end{equation*}
where $\bar{u}$ and $\upsilon$ are required to be at least $C^0$ continuous. The second term on the right-hand-side of the equation vanishes if Dirchlet boundary conditions are used.

The variational becomes
\begin{equation}
\label{eq:Variational}
I(u)  = \int_{x_{j-1/2}}^{x_{j+1/2}} \left ( \dfrac{h^3}{3} \dfrac{\partial \bar{u}}{\partial x} \dfrac{d\upsilon}{dx} + \bar{u}h\upsilon - G\upsilon \right ) \, dx = 0
\end{equation}
for an element

In the direct approximation of the finite element method, the desired function $q(x)$ is approximated by a weighted finite series
\begin{equation*}\label{eq:Hirsch-5.1.1}
q(x) \approx \hat{q}(x) = \sum\limits_{j=1}^N w_j(x) q_j.
\end{equation*}
involving the unknown nodal values, $q_j$ of the desired function. The set of locally defined piecewise functions, $w$ are also known as basis functions and $N$ is the total number of elements in the computational domain. In the Galerkin weighted residual method $\upsilon(x) = w(x)$.

\begin{figure}[h]
\centering
\includegraphics[width=12.0cm,keepaspectratio=true]{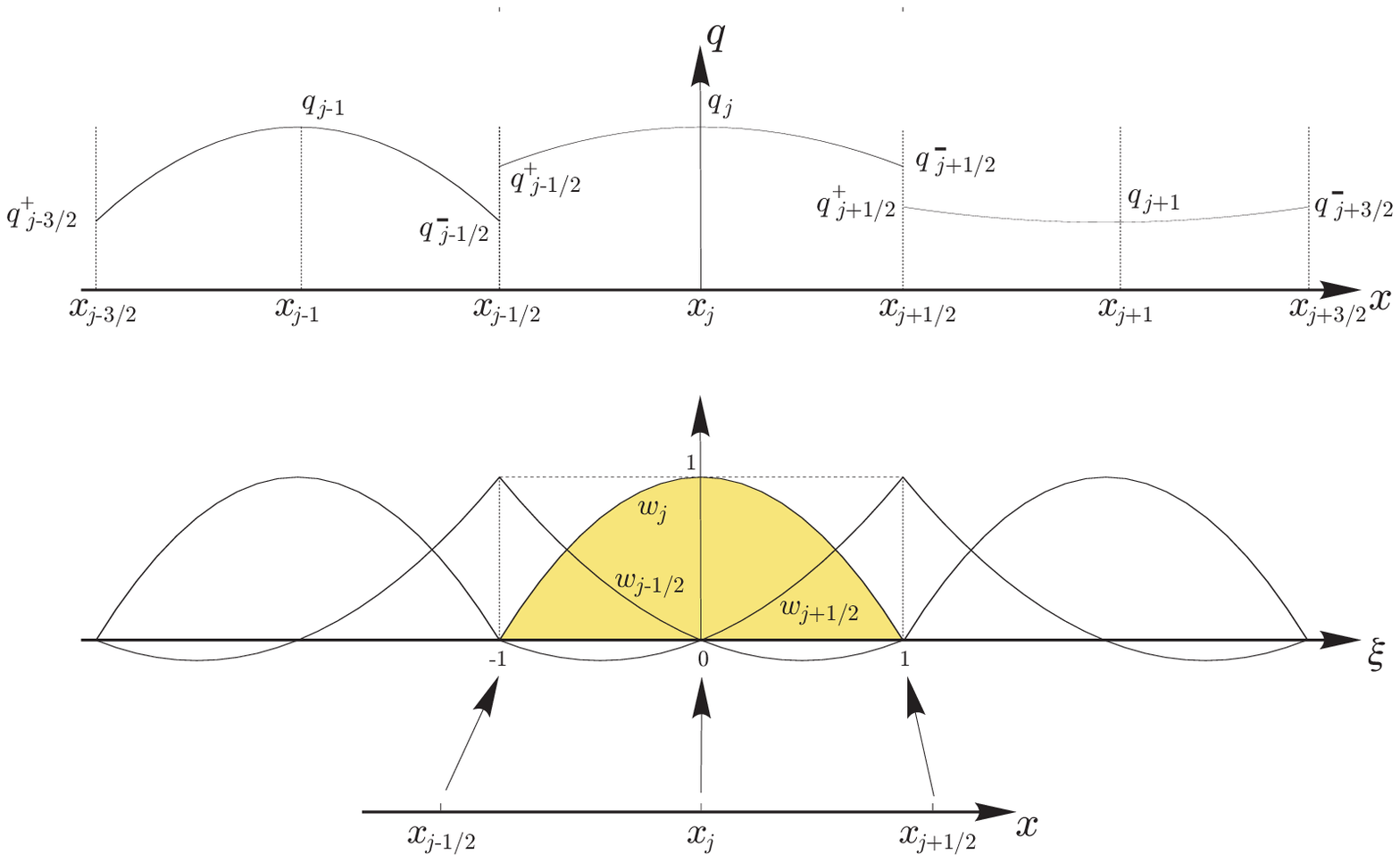}
\caption{The quadratic basis function used in the Galerkin finite elements which coencides with the finite volume cell.}
\label{fig:Quadratic_basis}
\end{figure}

For a third-order scheme, quadratic basis functions are used in each element, shown in Figure \ref{fig:Quadratic_basis}, spanning three nodes, $j-1/2$, $j$ and $j+1/2$. The finite elements
\begin{equation*}
\begin{array}{cccccc}
[x_{-1/2},x_{1/2}], & [x_{1/2},x_{3/2}], & \dots & [x_{j-1/2},x_{j+
+1/2}], & \dots & [x_{N-1/2},x_{N+1/2}] \\
I_1           &    I_2       &  \dots  & I_j  & \dots & I_N
\end{array}
\end{equation*}
coincides with the finite volume cells.

The Taylor series of a one-dimensional function, $q(\xi)$ and its derivatives including the coordinate transform mapping $[x_{j-1/2},x_j,x_{j+1/2}]$ to the $\xi$-space $[-1,0,1]$ is given in terms of the nodal values by
\begin{subequations}\label{eq:Iannelli-7.160}
\begin{equation*}
q_{j-1/2} = q(-1) = q(\xi) + (-1 - \xi)\dfrac{dq}{d\xi},
\end{equation*}
\begin{equation*}
q_j = q(0) = q(\xi) - \xi \dfrac{dq}{d\xi}
\end{equation*}
and
\begin{equation*}
q_{j+1/2} = q(1) = q(\xi) + (1 - \xi)\dfrac{dq}{d\xi}
\end{equation*}
\end{subequations}
where $x = \xi \Delta x/2 + x_j$ and $dx = d\xi \Delta x/2$ for a uniform grid. Solving for $q(\xi)$ and its derivative provides
\begin{equation*}
\left [
\begin{array}{ccc}
\dfrac{\xi(\xi - 1)}{2} & 1 - \xi^2 & \dfrac{\xi(\xi + 1)}{2} \\
\\
\dfrac{2\xi - 1}{2} & -2 \xi & \dfrac{2\xi + 1}{2}
\end{array}
\right ]
\left [
\begin{array}{c}
q_{j-1/2} \\
\\
q_j \\
\\
q_{j+1/2}
\end{array}
\right ] =
\left [
\begin{array}{c}
q \\
\\
\dfrac{dq}{d\xi}
\end{array}
\right ].
\end{equation*}
The coefficient matrix represents the weights for the quadratic interpolation of a quantity and its first derivative in space by the first and second rows respectively. The first row also represents the isometric mapping between $x$ and $\xi$.
\begin{equation*}
\left [
\begin{array}{ccc}
\dfrac{\xi(\xi - 1)}{2} & 1 - \xi^2 & \dfrac{\xi(\xi + 1)}{2} \\
\\
\dfrac{2\xi - 1}{2} & -2 \xi & \dfrac{2\xi + 1}{2}
\end{array}
\right ] =
\left [
\begin{array}{ccc}
w_{j-1/2} & w_j & w_{j+1/2} \\
\\
\dfrac{dw_{j-1/2}}{d\xi} & \dfrac{dw_j}{d\xi} & \dfrac{dw_{j+1/2}}{d\xi}
\end{array}
\right ] =
\left [
\begin{array}{c}
\mathbf{w}\\
\\
\dfrac{d\mathbf{w}}{d\xi}  \\
\\
\dfrac{d^2\mathbf{w}}{d\xi^2}
\end{array}
\right ].
\end{equation*}
For consistency, if $h_{x_{j-1}} = h_{x_j} = h_{x_{j+1}}$, then $w_{j+1/2} + w_j + w_{j+1/2} = 1$. Which is satisfied by the elements in the first row of the above matrix.

Consider each term in \eqref{eq:Variational},
\begin{gather*}
\int_{x_{j-1/2}}^{x_{j+1/2}} \dfrac{h^3}{3} \dfrac{\partial \bar{u}}{\partial x} \dfrac{\partial w}{\partial x} \, dx = \dfrac{dx}{d\xi} \int_{-1}^1 \left [ \dfrac{h(x_j + \xi)^3}{3} \dfrac{\partial \bar{u}(x_i + \xi)}{\partial \xi} \dfrac{d\xi}{dx} \dfrac{dw(x_j + \xi)}{d \xi} \dfrac{d\xi}{dx}  \right ] \, d\xi \\
= \dfrac{\Delta x}{2} \int_{-1}^1 \left ( \dfrac{1}{3} \left ( h^+_{j-1/2} w_{j-1/2} + h_jw_j + h^-_{j+1/2}w_{j+1/2} \right )^3 \right . \\
 \dfrac{2}{\Delta x} \left ( \bar{u}_{j-1/2} \dfrac{dw_{j-1/2}}{d \xi} + \bar{u}_j \dfrac{dw_j}{d \xi} + \bar{u}_{j+1/2} \dfrac{d w_{j+1/2}}{d \xi} \right )
 \dfrac{2}{\Delta x} \left .
 \left [
  \begin{array}{c}
    \dfrac{d w_{j-1/2}}{d \xi} \\
    \dfrac{d w_j}{d \xi} \\
    \dfrac{d w_{j+1/2}}{d \xi} \\
\end{array}
 \right ] \right ) \, d\xi \\
\end{gather*}
which can be written as $\mathbf{Q}_e \bar{\mathbf{u}}_e$, where the stiffness matrix, $\mathbf{Q}_e$ is given in the Appendix and $\bar{\mathbf{u}}_e = [ \bar{u}_{j-1/2}, \bar{u}_j, \bar{u}_{j+1/2} ]^T$.

Consider the term
\begin{gather*}
\int_{x_{j-1/2}}^{x_{j+1/2}} \bar{u} h w \, dx = \dfrac{dx}{d\xi} \int_{-1}^1 \left ( \bar{u}(x_i + \xi) h(x_j + \xi) w(x_j + \xi)  \right ) \, d\xi \\
= \dfrac{\Delta x}{2} \int_{-1}^1  \left ( \bar{u}_{j-1/2}w_{j-1/2} + \bar{u}_j w_j + \bar{u}_{j+1/2} w_{j+1/2} \right )  \left ( h^+_{j-1/2} w_{j-1/2} + h_jw_j + h^-_{j+1/2}w_{j+1/2} \right )
 \left[
  \begin{array}{c}
    w_{j-1/2} \\
    w_j \\
    w_{j+1/2} \\
  \end{array}
\right ]
 \,\, d\xi. \\
\end{gather*}
Which can be written as $\mathbf{P}_e \bar{\mathbf{u}}_e$ where
\begin{equation*}
\mathbf{P}_e  = \dfrac{\Delta x}{420} \left [
  \begin{array}{ccc}
   39 h^+_{j-1/2}  +  20  h_j   - 3 h^-_{j+1/2} &
   20 h^+_{j-1/2}  + 16  h_j  - 8 h^-_{j+1/2} &
   -3 h^+_{j-1/2} - 8  h_j  - 3 h^-_{j+1/2} \\
    20 h^+_{j-1/2} +  16 h_j  - 8 h^-_{j+1/2} &
   16 h^+_{j-1/2} + 192 h_j + 16 h^-_{j+1/2} &
   -8  h^+_{j-1/2} + 16 h_j + 20 h^-_{j+1/2} \\
      -3 h^+_{j-1/2} - 8 h_j  - 3 h^-_{j+1/2} &
   -8 h^+_{j-1/2} + 16 h_j + 20 h^-_{j+1/2} &
   -3 h^+_{j-1/2} + 20 h_j + 39 h^-_{j+1/2} \\
\end{array}
 \right ].
\end{equation*}

The final term
\begin{gather*}
\int_{x_{j-1/2}}^{x_{j+1/2}} G(x) w(x) \,\, dx = \dfrac{dx}{d\xi} \int_{-1}^1 G(x_j + \xi)
\left[
  \begin{array}{c}
    w_{j-1/2} \\
    w_j \\
    w_{j+1/2} \\
  \end{array}
\right ]
 \,\, d\xi \\
= \dfrac{\Delta x}{2} \int_{-1}^1 \left ( G_{j-1/2}^+(\xi) w_{j-1/2}(\xi) + G_j(\xi) w_j(\xi) + G_{j+1/2}^-(\xi) w_{j+1/2}(\xi) \right )
\left[
  \begin{array}{c}
    w_{j-1/2} \\
    w_j \\
    w_{j+1/2} \\
  \end{array}
\right ]
 \,\, d\xi
\end{gather*}
Evaluating the integrals yields
\begin{align*}
\mathbf{R}_e = \left[
                 \begin{array}{c}
                   r_1 \\
                   r_2 \\
                   r_3 \\
                 \end{array}
               \right] = \dfrac{\Delta x}{30} \left[
                                                \begin{array}{c}
                                                  4  G_{j-1/2}^+ + 2 G_j  - G_{j+1/2}^-\\
                                                  2  G_{j-1/2}^+ + 16 G_j + 2 G_{j+1/2}^-\\
                                                  -G_{j-1/2}^+ + 2 G_j + 4 G_{j+1/2}^-\\
                                                \end{array}
                                                \right ]
\end{align*}

For each element
\begin{equation*}
\left [ \mathbf{Q}_e + \mathbf{P}_e \right ] \bar{\mathbf{u}}_e -\mathbf{R}_e = 0
\end{equation*}
The assembled system of equations can be written in matrix form,
\begin{equation}
\label{eq:matrix_FV}
[\mathbf{A}][\mathbf{u}] = [\mathbf{R}]
\end{equation}
where the stiffness matrix $\mathbf{A} = \sum\limits_e [\mathbf{Q}_e + \mathbf{P}_e]$,  the load matrix $\mathbf{R} = \sum\limits_e [\mathbf{R}_e]$ and the vector of unknowns $\mathbf{u} = \sum\limits_e [\mathbf{u}_e]$ with $\mathbf{u}_e^T = [\bar{u}_{j-1/2} \quad \bar{u}_j \quad \bar{u}_{j+1/2}]$.

Written for all computational nodes, the coefficient matrix $[\mathbf{A}]$ is a penta-diagonal matrix, see the Appendix and the system of equations, \eqref{eq:matrix_FV} is straightforward to solve for the primitive variable, $\mathbf{u}_e$,

The conservative quantities, $h$ and $G$ can be discontinuous across each element, see Figure \ref{fig:Quadratic_basis}, whilst the unknown quantity, $\bar{u}$ is assumed to be continuous and approximated by a piecewise quadratic function. The values $\bar{u}_j$ and $\bar{u}_{j \pm 1/2}$ are use to calculate the gradient used in the second component in \eqref{eq:Serre_flux}, shown in \eqref{eq:local_flux_discrete} without a loose in accuracy.

\subsection{Ensembles Scheme}

The third-order three-stage strong stability preserving Runge-Kutta scheme solution of the Serre equations involves the following steps
\begin{gather}
\underbrace{\bm{\bar{q}}^n_j \stackrel{\mathcal{H}}{\longrightarrow} \bar{\bm{u}}^n_j}_{\text{\ding{192}}}  \rightarrow \underbrace{\bm{\bar{q}}^{(1)}_j
        =\bm{\bar{q}}^{n}_j - \dfrac{\Delta t}{\Delta x} \left ( \bm{F}_{j+1/2}^n - \bm{F}_{j-1/2}^n \right)}_{\text{\ding{193} First Euler Step}} \nonumber \\
 \underbrace{ \color{black}\bm{\bar{q}}^{(1)}_j \stackrel{\mathcal{H}}{\longrightarrow} \bar{\bm{u}}^{(1)}_j}_{\text{ \ding{194}}} \rightarrow
        \underbrace{ \bm{\bar{q}}^{(2)}_j
         = \bm{\bar{q}}^{(1)}_j - \dfrac{\Delta t}{\Delta x} \left ( \bm{F}_{j+1/2}^{(1)} - \bm{F}_{j-1/2}^{(1)} \right)}_{\text{\ding{195} Second Euler Step}} \nonumber \\
 \color{black} \underbrace{\bm{\bar{q}}^{(3)}_j
         = \dfrac{3}{4}\bm{\bar{q}}^n_j + \dfrac{1}{4}\bm{\bar{q}}^{(2)}_j}_{\text{\ding{196} Intermediate Step}}\nonumber \\
         \underbrace{\bm{\bar{q}}^{(3)}_j \stackrel{\mathcal{H}}{\longrightarrow} \bar{\bm{u}}^{(3)}_j}_{\text{\ding{197}}}  \rightarrow \underbrace{\bm{\bar{q}}^{(4)}_j
        =\bm{\bar{q}}^{(3)}_j - \dfrac{\Delta t}{\Delta x} \left ( \bm{F}_{j+1/2}^{(3)} - \bm{F}_{j-1/2}^{(3)} \right)}_{\text{\ding{198} Third Euler Step}} \nonumber \\
  \color{black} \underbrace{\mathbf{\bm{\bar{q}}}^{n+1}_j
         = \dfrac{1}{3}\bm{\bar{q}}^n_j + \dfrac{2}{3}\bm{\bar{q}}^{(4)}_j}_{\text{\ding{199} Averaging Step}}\nonumber
\end{gather}

\begin{description}
\item [Step 1:] Given $\bm{\bar{q}_j} = [\bar{h}_j, \bar{G}_j]$ the remaining primitive variable $\bar{\bm{u}} = \mathcal{H}[\bm{\bar{q}}]$  is obtained by solving the second-order elliptic equation, \eqref{eq:G_FE} using the third-order finite element scheme. The operator $\mathcal{H}$ involves the following steps
\begin{gather*}
\mathcal{H}(\bm{q}) = \left \lbrace \begin{matrix} \bm{\bar{q}} \stackrel{\mathcal{M}}{\rightarrow} \bm{q}_j \\
\bar{\bm{q}} \stackrel{\mathcal{A}}{\rightarrow} (\bm{q}_{j+1/2}) \\
\bm{q}_{j+1/2} \stackrel{\mathcal{F}}{\rightarrow} \mathbf{u}_e
\end{matrix}
\right .
\end{gather*}
where $\mathcal{A}$ is the reconstruction of the point values at the cell interface using the Koren limiter and $\mathcal{F}$ is the third-order finite element solution of the second-order elliptic equation for the point values $\mathbf{u}_e$ using \eqref{eq:matrix_FV}.

\item[Step 2:] Perform the reconstruction using the cell averages, $\bar{\bm{q}} = [\bar{h}_j, \bar{G}_j]^T$ and solve the local Riemann problem to obtain the flux $\bm{F}_{j \pm 1/2}$ of material across a cell interface. Evolve the solution using a first-order Euler time integration for the conserved quantities, $\bar{\bm{q}}_j$.

\item[Steps 3 and 4:] Repeat Steps $1$ and $2$ with the values obtained from Step $2$ and evolve using another first-order Euler step.

\item[Step 5:] Calculate weighted intermediate values using the initial values and the results from Step $4$.

\item[Steps 6 and 7:] Repeat Steps $1$ and $2$ using the intermediate values.

\item[Step 8:]
The solution at the next time level is obtained by averaging the initial values and the values obtained from the third Euler step. This completes the third-order strong stability preserving Runge-Kutta time integration.
\end{description}

\section{Numerical Simulations}

Convergence rate of the proposed schemes is determined using a known analytical solution to the Serre equations. Data from two laboratory experiments are used to validate the proposed models and the simulation of the dam-break problem is used to show that the model is stable for simulating a wide range of flow problems.

\subsection{Analytical Solutions}

Simulating the propagation of solitons is a common test for Boussinesq-type equations. The Serre equations, \eqref{eq:Serre_conservative_form} has the following analytical  soliton solution\cite{El-etal-2006}
\begin{subequations}\label{eq:Carter-Cienfuegos-solitary-wave}
\label{eq:Carter-9}
\begin{equation}
h(x,t) = a_0 + a_1 \text{sech}^2(\kappa(x - ct))
\end{equation}
and
\begin{equation}
\bar{u}(x,t) = c \left ( \dfrac{h(x,t) - a_0}{h(x,t)} \right )
\end{equation}
\end{subequations}
with $\kappa = \sqrt{3a_1}/({2a_0\sqrt{a_0 + a_1}})$ and $c = \sqrt{g(a_0 + a_1)}$.

Solitary waves propagate at constant speed without deformation. Therefore, there is a balance between nonlinear and dispersive effects.  A numerical scheme must accurately model the equilibrium between amplitude and frequency dispersion in order to simulate the  propagation of the wave profile at constant shape and speed. A result of a poorly balanced numerical scheme is the simulation of trailing edge dispersion waves which cause a reduction in wave height and celerity.

The results from a numerical scheme are compared to the corresponding analytical solution by using the non-dimensionless $L_1$ norm
\begin{equation*}
L_1(q_j,q(x_j)) = \dfrac{\sum_{j=1}^m |q_j - q(x_j)|}{\sum_{j=1}^m |q(x_j)|}
\label{eq:Li_norm}
\end{equation*}
where, $q_j$ is the simulated values of $q(x,t)$ at $x_j$, and $q(x_j)$ is the corresponding analytical solution. The $L_1$ norm is calculated using all the computational nodes, $j = 1,\ldots, m$.

The prototypical example is a solitary wave predicted by \eqref{eq:Carter-Cienfuegos-solitary-wave} with, $a_0 = 10$m, an amplitude of $a_1 = 1.0$m. The soliton has a celerity, $c = 10.387974$m/s and $\kappa = 0.026112$/m.

\begin{figure}[ht]
\centering
\begin{tabular}{ccc}
\begin{psfrags}%
\psfragscanon%
%
\psfrag{s03}[t][t]{\setlength{\tabcolsep}{0pt}\begin{tabular}{c}{\Large $x$(m)}\end{tabular}}%
\psfrag{s04}[b][b]{\setlength{\tabcolsep}{0pt}\begin{tabular}{c}{\Large $h$(m)}\end{tabular}}%
%
\psfrag{x01}[t][t]{-500}%
\psfrag{x02}[t][t]{0}%
\psfrag{x03}[t][t]{500}%
\psfrag{x04}[t][t]{1000}%
\psfrag{x05}[t][t]{1500}%
%
\psfrag{v01}[r][r]{9.8}%
\psfrag{v02}[r][r]{10.0}%
\psfrag{v03}[r][r]{10.2}%
\psfrag{v04}[r][r]{10.4}%
\psfrag{v05}[r][r]{10.6}%
\psfrag{v06}[r][r]{10.8}%
\psfrag{v07}[r][r]{11.0}%
\psfrag{v08}[r][r]{11.2}%
%
\resizebox{6cm}{!}{\includegraphics{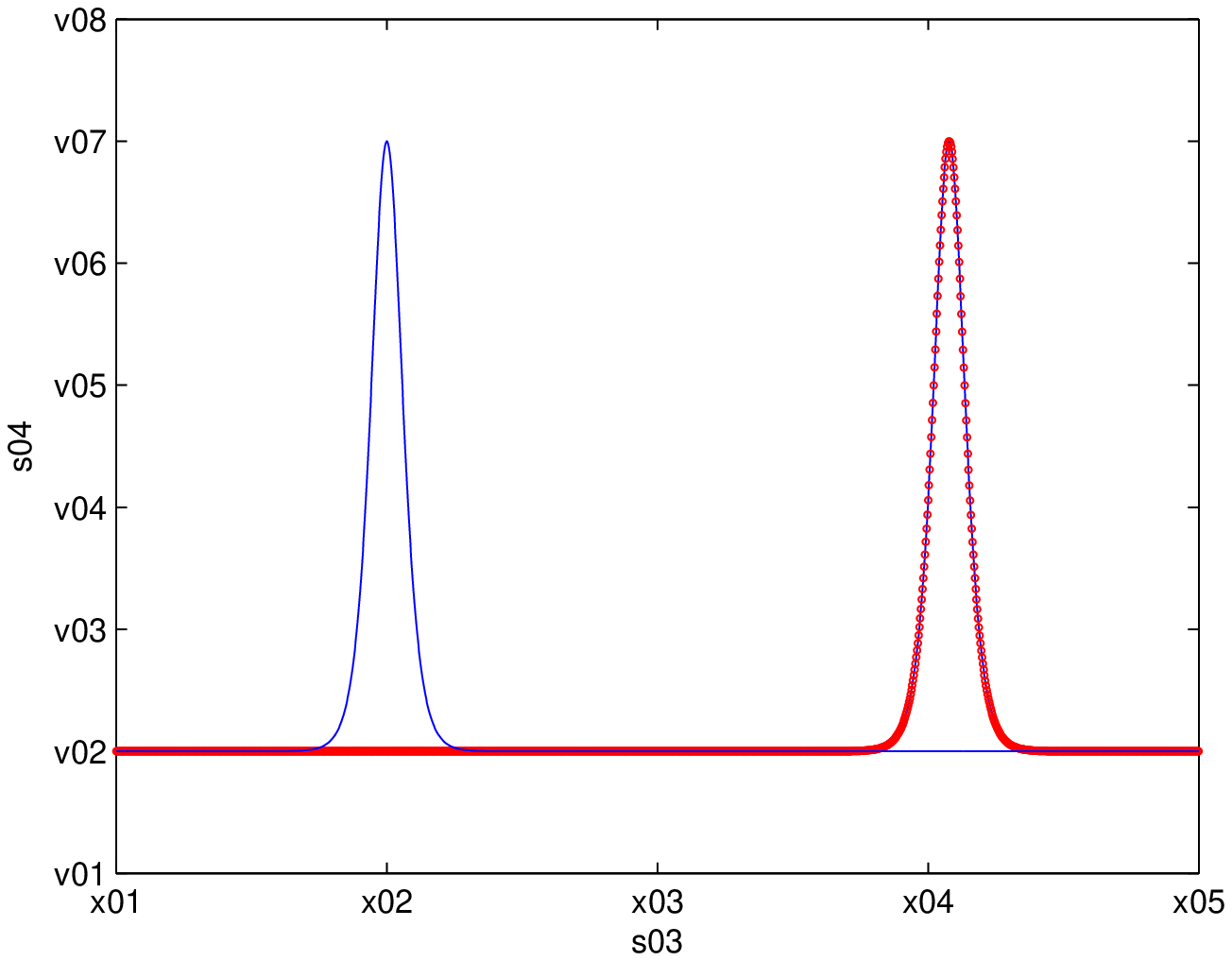}}%
\end{psfrags}%
 & & %
 \begin{psfrags}%
\psfragscanon%
%
\psfrag{s03}[t][t]{\setlength{\tabcolsep}{0pt}\begin{tabular}{c}{\Large $x$(m)}\end{tabular}}%
\psfrag{s04}[b][b]{\setlength{\tabcolsep}{0pt}\begin{tabular}{c}{\Large $\bar{u}$(m/s)}\end{tabular}}%
%
\psfrag{x01}[t][t]{-500}%
\psfrag{x02}[t][t]{0}%
\psfrag{x03}[t][t]{500}%
\psfrag{x04}[t][t]{1000}%
\psfrag{x05}[t][t]{1500}%
%
\psfrag{v01}[r][r]{-0.2}%
\psfrag{v02}[r][r]{0.0}%
\psfrag{v03}[r][r]{0.2}%
\psfrag{v04}[r][r]{0.4}%
\psfrag{v05}[r][r]{0.6}%
\psfrag{v06}[r][r]{0.8}%
\psfrag{v07}[r][r]{1.0}%
%
\resizebox{6cm}{!}{\includegraphics{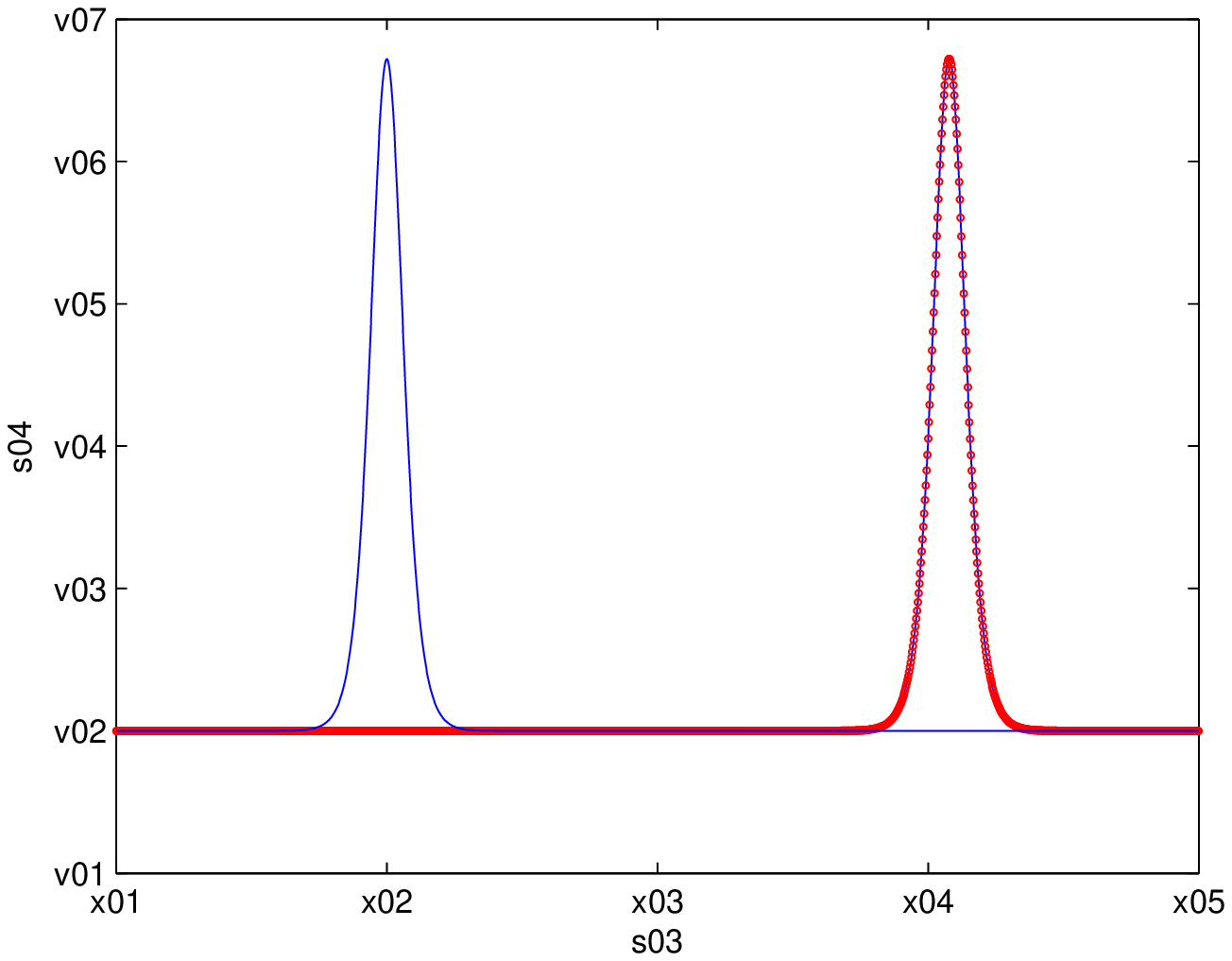}}%
\end{psfrags}%
\\
($a$) & & ($b$)
\end{tabular}
\caption{The progress of a initial solitary wave,  given by \eqref{eq:Carter-Cienfuegos-solitary-wave} over a horizontal bed predicted by the third-order scheme solution of \eqref{eq:Serre_Vector} ($\circ$), where $\Delta x = ??????$m, and $Cr = 0.1$, at $t = 100$s with the water depth, $h(x,t)$ shown in ($a$) and the velocity, $\bar{u}(x,t)$ in ($b$) plotted against the analytical solution (\textemdash).}
\label{fig:Solitary_wave_simulation_second-order}
\end{figure}

The boundary conditions imposed on the models are $a_0(x)= 10$m and $\bar{u}_0(x) = 0$m/s at the upstream and downstream boundaries. The model parameters are; $\Delta x = 0.1$m, $Cr = 0.1$ and $\Delta t = Cr/\Delta x$s. Using these parameters, the initial soliton profile and velocity, the analytical and the simulated water depth and velocity at $t = 100$s are shown in Figure \ref{fig:Solitary_wave_simulation_second-order} for the third-order scheme.  The soliton amplitude is accurately predicted and the soliton speed is captured correctly.

\begin{figure}[htb]
\centering
\begin{psfrags}%
\psfragscanon%
%
\psfrag{s02}[l][l]{\setlength{\tabcolsep}{0pt}\begin{tabular}{l}$1$\end{tabular}}%
\psfrag{s03}[l][l]{\setlength{\tabcolsep}{0pt}\begin{tabular}{l}$3$\end{tabular}}%
\psfrag{s04}[r][r]{}%
\psfrag{s05}[t][t]{\setlength{\tabcolsep}{0pt}\begin{tabular}{c}{\Large $Log_{10}\Delta x$}\end{tabular}}%
\psfrag{s06}[b][b]{\setlength{\tabcolsep}{0pt}\begin{tabular}{c}{\Large $Log_{10}L_1$}\end{tabular}}%
%
\psfrag{x01}[t][t]{$10^{-3}$}%
\psfrag{x02}[t][t]{$10^{-2}$}%
\psfrag{x03}[t][t]{$10^{-1}$}%
\psfrag{x04}[t][t]{$10^{0}$}%
\psfrag{x05}[t][t]{$10^{1}$}%
\psfrag{x06}[t][t]{$10^{2}$}%
\psfrag{x07}[t][t]{$10^{3}$}%
%
\psfrag{v01}[r][r]{$10^{-12}$}%
\psfrag{v02}[r][r]{$10^{-10}$}%
\psfrag{v03}[r][r]{$10^{-8}$}%
\psfrag{v04}[r][r]{$10^{-6}$}%
\psfrag{v05}[r][r]{$10^{-4}$}%
\psfrag{v06}[r][r]{$10^{-2}$}%
\psfrag{v07}[r][r]{$10^{0}$}%
\psfrag{v08}[r][r]{$10^{2}$}%
%
\resizebox{6cm}{!}{\includegraphics{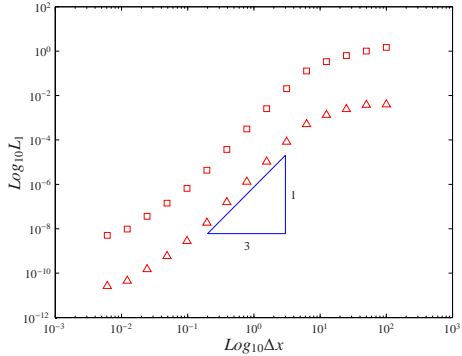}}%
\end{psfrags}%
\caption{The $L_1$ convergence rate for the simulated water depth ($\triangle$) and velocity ($\square$) obtained from the second-order finite volume central upwind scheme solution of \eqref{eq:Serre_Vector} to the solitary wave example,  given by \eqref{eq:Carter-Cienfuegos-solitary-wave}.}
\label{fig:Convergence_Soliton_second-order}
\end{figure}

Performing the simulation for a range of $\Delta x$ and keeping $Cr = 0.1$, the $L_1$ norm between the simulated and analytical solution was calculated for the water depth and fluid velocity. Plotting the $\log_{10}L_1$ against $\log_{10}\Delta x$ reveals that the proposed strategy for solving the Serre equations is third-order accurate, see Figure \ref{fig:Convergence_Soliton_second-order}. Here, the usual convergence behaviour is observed. For coarse grids, the numerical scheme is not able to resolve the oscillations and the true convergence rate of the scheme is not realized. As the grid is refined convergence is confirmed and for very small grid spacing the error deteriorates because it is dominated by the accumulation of roundoff errors.

Clearly, for the simulation of the smooth soliton problem, the third-order schemes is capable of predicting the soliton speed and its amplitude.

\subsection{Labratory Experiments}\label{Laboratory_Experiments}

The frictionless horizontal flume experiment from Hammack and Segur\cite{Hammack-Segur-1978-337}, involving a negative amplitude rectangular wave and the more recent surge propagation experiment conducted by Chanson\cite{Chanson-H-2009-104}, are used to validate the hird-order schemes described in Section \ref{numerical_scheme}. Both produce highly dispersive waves from an abrupt change in the initial flow conditions. In these experiments the non-hydrostatic terms cannot be neglected in the momentum equation.

\subsubsection{Undular Bore}\label{Undular Bore}

An undular bore was created in a  large tilting flume at the Civil Engineering Department, University of Queensland. The channel is $0.5$m wide, $12$m in length and the undular bore was created in the horizontal flume, which has a smooth PVC bed and glass walls. A radial gate located at the downstream end of the flume, $x = 11.9$m controls the water depth in the flume. The radial gate is used during the experiments to produce steady subcritical flow in the flume which remains constant for the duration of the experiment. Steady flow condition are established for 15 minutes prior to an experiment. Adjacent to the radial gate is a rapidly closing Tainter gate at, $x = 11.15$m that spans the full width of the flume.  An undular bore is generated by the rapid closure of the Tainter gate, which is estimated to take less than $0.2$s, when water accumulates at the Tainer gate forming an upstream progressing undular bore. The experiment ceases when the bore reaches the upstream intake structure to avoid any interference from wave reflection. Acoustic displacement meters, located at the flume centerline at; $x = 8.0$, $6.0$, $5.0$, $4.55$, $4.0$ and $3.0$m record the progress of the bore and dispersive waves with time. Data acquisition starts 30 seconds prior to the closure of the Tainter gate.

The boundary conditions imposed in all the models are; at the upstream boundary,  $h(0,t) = 0.192$m and $\bar{u}(0,t) = 0.199$m$^3$/s and at the downstream Tainer gate, $h(11.15,t) = 0.22$m and $\bar{u}(11.5,t) = 0$m/s. In all the simulations, $\Delta x = 0.01115$m and $Cr = 0.2$ in the third-order scheme.

\begin{figure}[htb]
\centering
\begin{tabular}{cc}
\begin{psfrags}%
\psfragscanon%
%
\psfrag{s03}[t][t]{\setlength{\tabcolsep}{0pt}\begin{tabular}{c}{\Large $t$(s)}\end{tabular}}%
\psfrag{s04}[b][b]{\setlength{\tabcolsep}{0pt}\begin{tabular}{c}{\Large $h$(m)}\end{tabular}}%
%
\psfrag{x01}[t][t]{25}%
\psfrag{x02}[t][t]{30}%
\psfrag{x03}[t][t]{35}%
\psfrag{x04}[t][t]{40}%
%
\psfrag{v01}[r][r]{0.18}%
\psfrag{v02}[r][r]{0.19}%
\psfrag{v03}[r][r]{0.20}%
\psfrag{v04}[r][r]{0.21}%
\psfrag{v05}[r][r]{0.22}%
\psfrag{v06}[r][r]{0.23}%
\psfrag{v07}[r][r]{0.24}%
%
\resizebox{6cm}{!}{\includegraphics{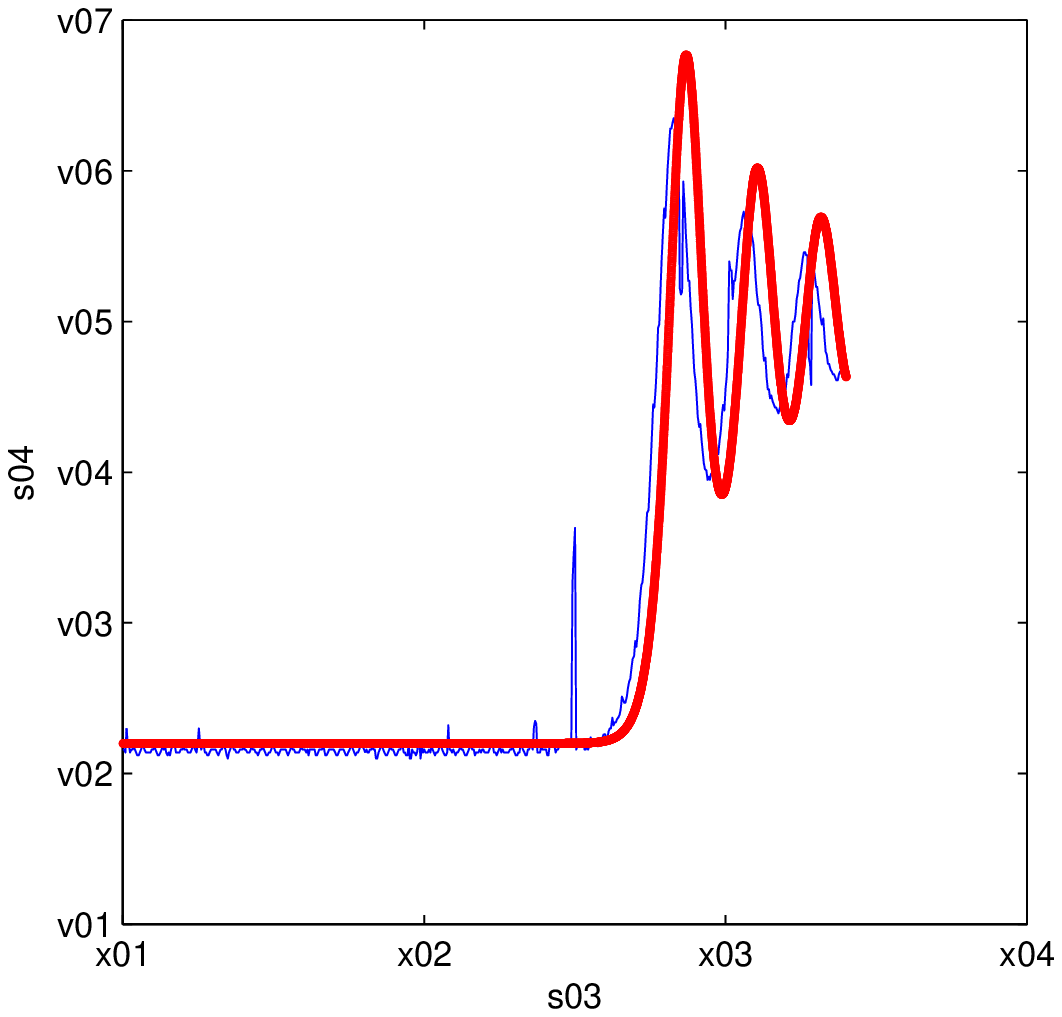}}%
\end{psfrags}%
& %
\begin{psfrags}%
\psfragscanon%
%
\psfrag{s03}[t][t]{\setlength{\tabcolsep}{0pt}\begin{tabular}{c}{\Large $t$(s)}\end{tabular}}%
\psfrag{s04}[b][b]{\setlength{\tabcolsep}{0pt}\begin{tabular}{c}{\Large $h$(m)}\end{tabular}}%
%
\psfrag{x01}[t][t]{25}%
\psfrag{x02}[t][t]{30}%
\psfrag{x03}[t][t]{35}%
\psfrag{x04}[t][t]{40}%
%
\psfrag{v01}[r][r]{0.18}%
\psfrag{v02}[r][r]{0.19}%
\psfrag{v03}[r][r]{0.20}%
\psfrag{v04}[r][r]{0.21}%
\psfrag{v05}[r][r]{0.22}%
\psfrag{v06}[r][r]{0.23}%
\psfrag{v07}[r][r]{0.24}%
%
\resizebox{6cm}{!}{\includegraphics{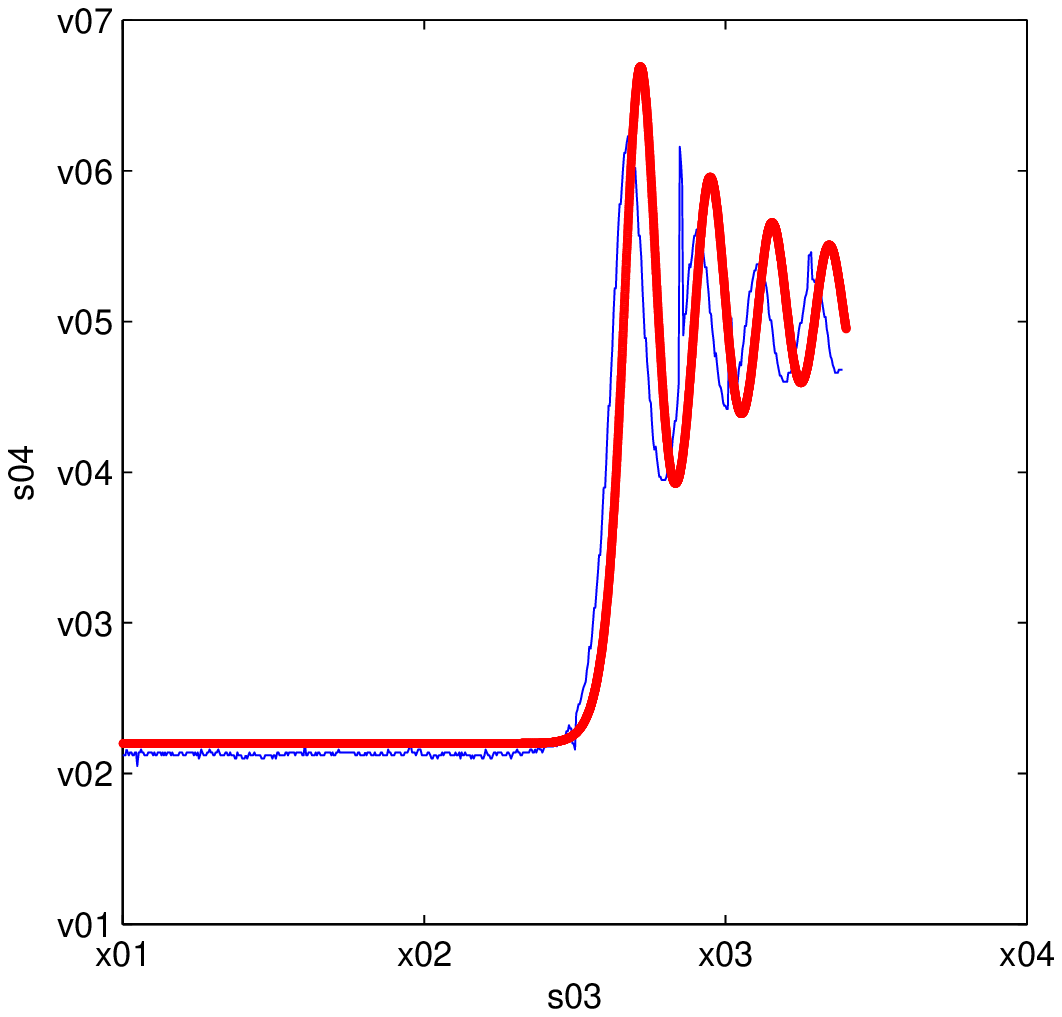}}%
\end{psfrags}%
\\
\phantom{x} & \\
($a$) & ($b$) \\
\phantom{x} & \\
\begin{psfrags}%
\psfragscanon%
%
\psfrag{s03}[t][t]{\setlength{\tabcolsep}{0pt}\begin{tabular}{c}{\Large $t$(s)}\end{tabular}}%
\psfrag{s04}[b][b]{\setlength{\tabcolsep}{0pt}\begin{tabular}{c}{\Large $h$(m)}\end{tabular}}%
%
\psfrag{x01}[t][t]{25}%
\psfrag{x02}[t][t]{30}%
\psfrag{x03}[t][t]{35}%
\psfrag{x04}[t][t]{40}%
%
\psfrag{v01}[r][r]{0.18}%
\psfrag{v02}[r][r]{0.19}%
\psfrag{v03}[r][r]{0.20}%
\psfrag{v04}[r][r]{0.21}%
\psfrag{v05}[r][r]{0.22}%
\psfrag{v06}[r][r]{0.23}%
\psfrag{v07}[r][r]{0.24}%
%
\resizebox{6cm}{!}{\includegraphics{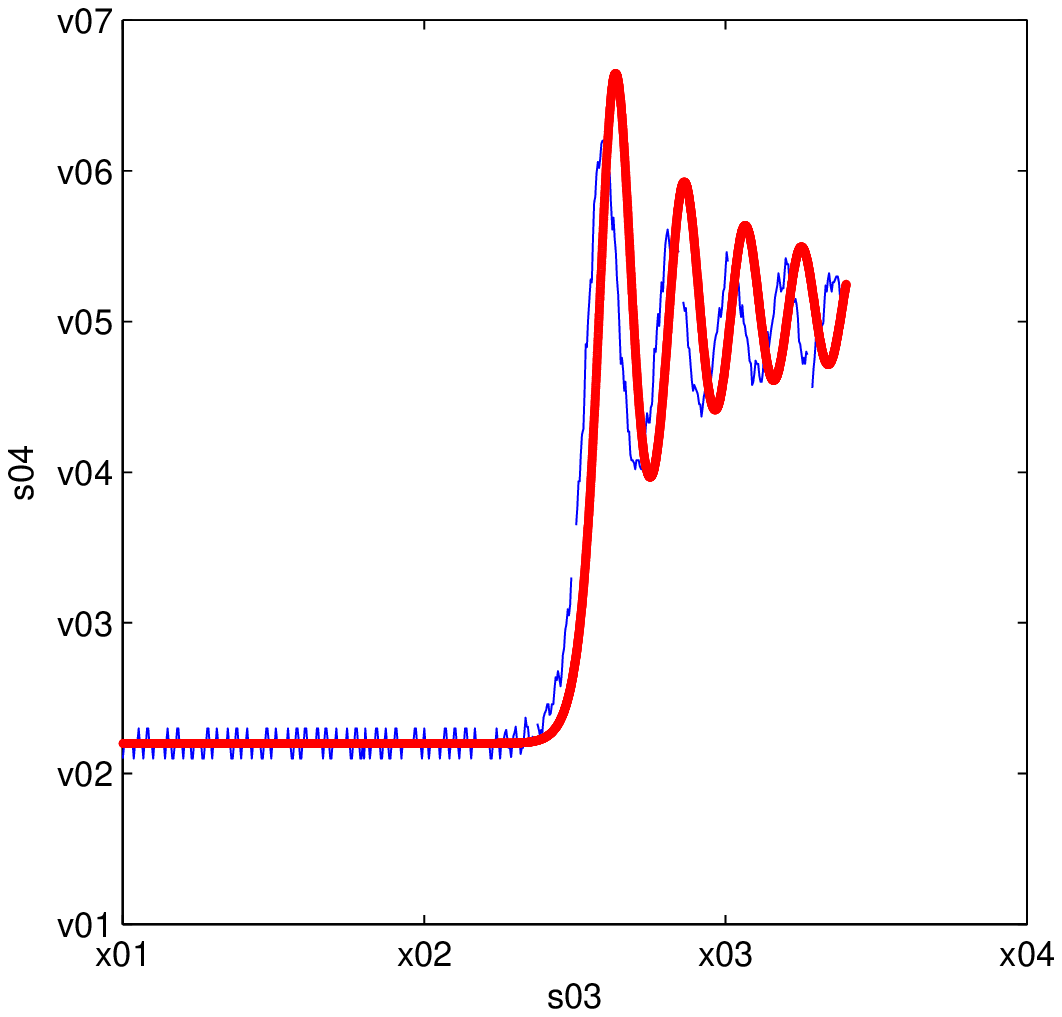}}%
\end{psfrags}%
&
\begin{psfrags}%
\psfragscanon%
%
\psfrag{s03}[t][t]{\setlength{\tabcolsep}{0pt}\begin{tabular}{c}{\Large $t$(s)}\end{tabular}}%
\psfrag{s04}[b][b]{\setlength{\tabcolsep}{0pt}\begin{tabular}{c}{\Large $h$(m)}\end{tabular}}%
%
\psfrag{x01}[t][t]{25}%
\psfrag{x02}[t][t]{30}%
\psfrag{x03}[t][t]{35}%
\psfrag{x04}[t][t]{40}%
%
\psfrag{v01}[r][r]{0.18}%
\psfrag{v02}[r][r]{0.19}%
\psfrag{v03}[r][r]{0.20}%
\psfrag{v04}[r][r]{0.21}%
\psfrag{v05}[r][r]{0.22}%
\psfrag{v06}[r][r]{0.23}%
\psfrag{v07}[r][r]{0.24}%
%
\resizebox{6cm}{!}{\includegraphics{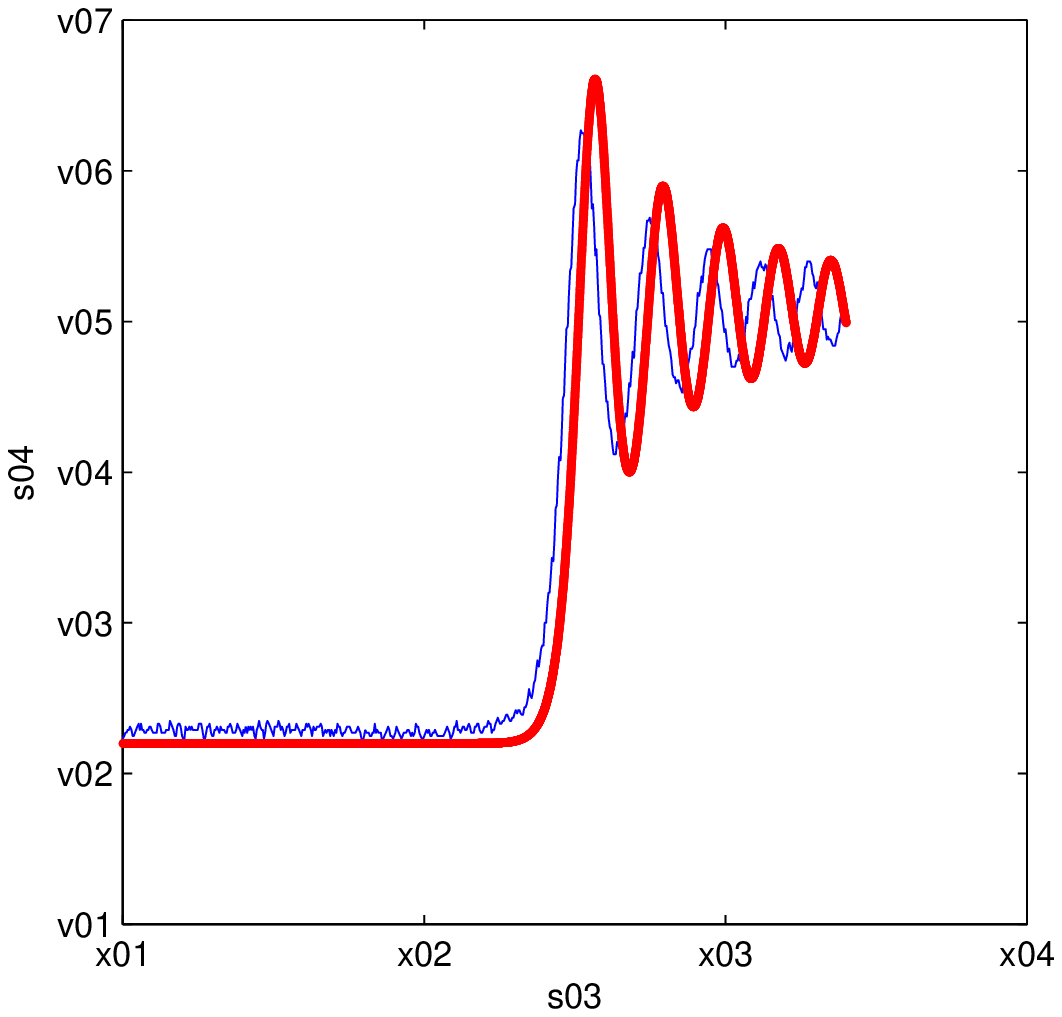}}%
\end{psfrags}%
\\
 \phantom{x} & \\
($c$) & ($d$) \\
\phantom{x} & \\
\begin{psfrags}%
\psfragscanon%
%
\psfrag{s03}[t][t]{\setlength{\tabcolsep}{0pt}\begin{tabular}{c}{\Large $t$(s)}\end{tabular}}%
\psfrag{s04}[b][b]{\setlength{\tabcolsep}{0pt}\begin{tabular}{c}{\Large $h$(m)}\end{tabular}}%
%
\psfrag{x01}[t][t]{25}%
\psfrag{x02}[t][t]{30}%
\psfrag{x03}[t][t]{35}%
\psfrag{x04}[t][t]{40}%
%
\psfrag{v01}[r][r]{0.18}%
\psfrag{v02}[r][r]{0.19}%
\psfrag{v03}[r][r]{0.20}%
\psfrag{v04}[r][r]{0.21}%
\psfrag{v05}[r][r]{0.22}%
\psfrag{v06}[r][r]{0.23}%
\psfrag{v07}[r][r]{0.24}%
%
\resizebox{6cm}{!}{\includegraphics{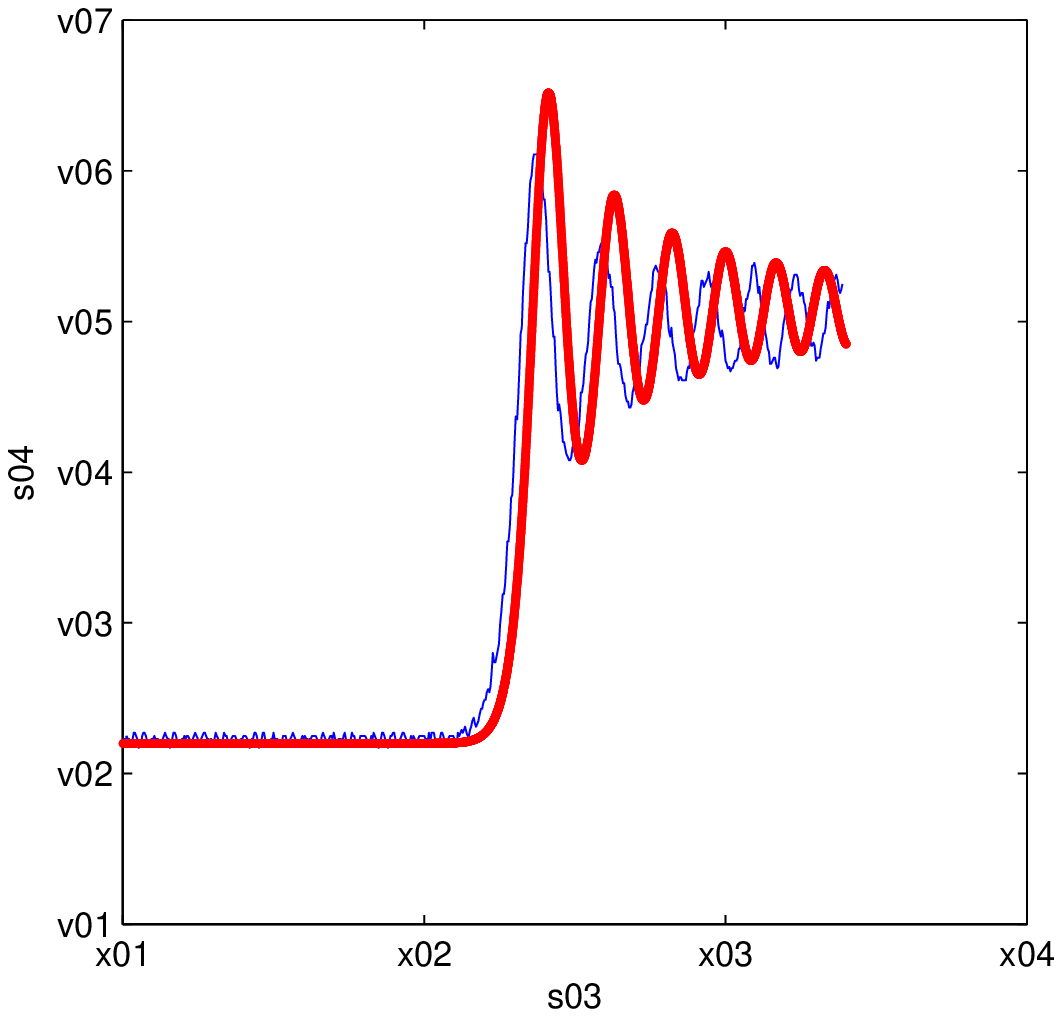}}%
\end{psfrags}%
&
\begin{psfrags}%
\psfragscanon%
%
\psfrag{s03}[t][t]{\setlength{\tabcolsep}{0pt}\begin{tabular}{c}{\Large $t$(s)}\end{tabular}}%
\psfrag{s04}[b][b]{\setlength{\tabcolsep}{0pt}\begin{tabular}{c}{\Large $h$(m)}\end{tabular}}%
%
\psfrag{x01}[t][t]{25}%
\psfrag{x02}[t][t]{30}%
\psfrag{x03}[t][t]{35}%
\psfrag{x04}[t][t]{40}%
%
\psfrag{v01}[r][r]{0.18}%
\psfrag{v02}[r][r]{0.19}%
\psfrag{v03}[r][r]{0.20}%
\psfrag{v04}[r][r]{0.21}%
\psfrag{v05}[r][r]{0.22}%
\psfrag{v06}[r][r]{0.23}%
\psfrag{v07}[r][r]{0.24}%
%
\resizebox{6cm}{!}{\includegraphics{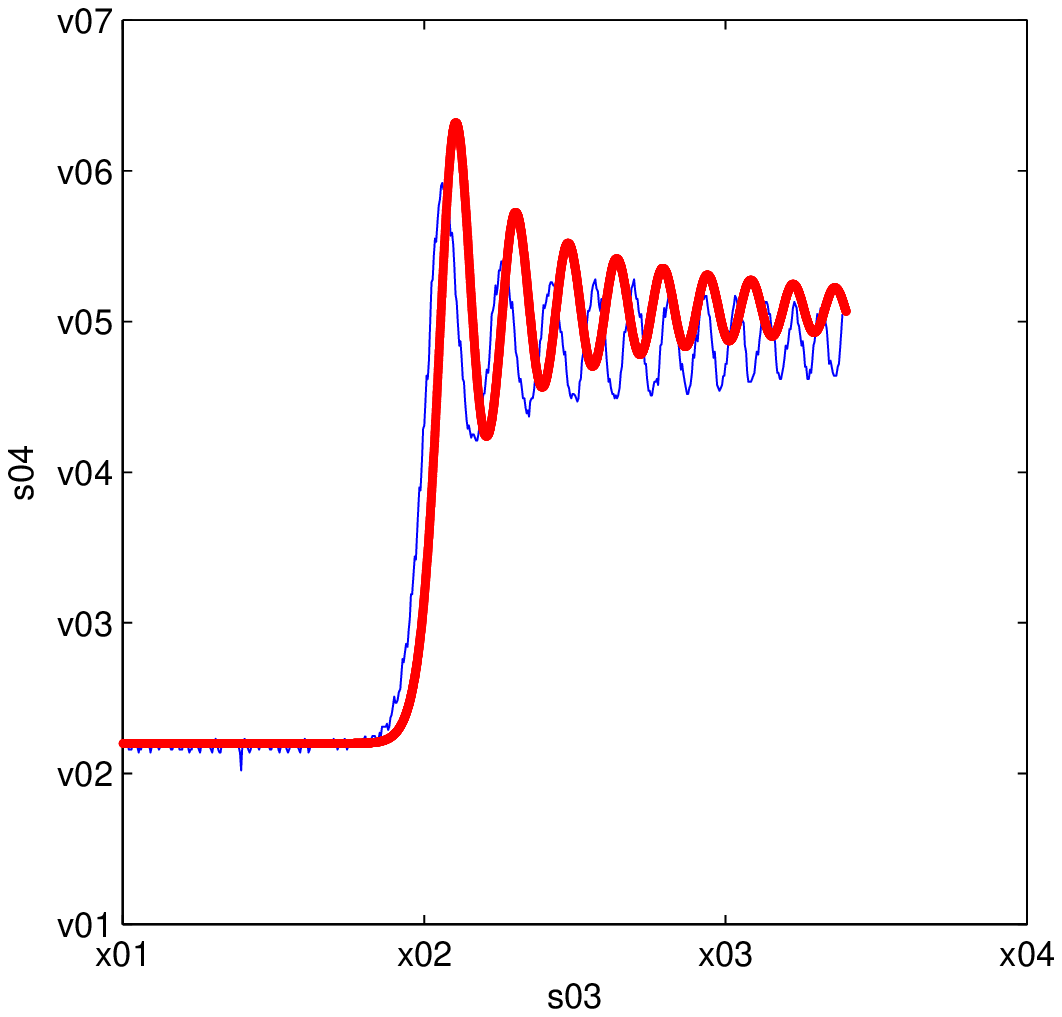}}%
\end{psfrags}  \\%
\phantom{x} & \\
($e$)  & ($f$) \\
\phantom{x}  \\
\end{tabular}
\caption{Measured (\textemdash) and simulated ($\circ$) water depth, $h(x,t)$ for the undular bore experiment in a frictionless rectangular channel using the second-order central and Euler time integration solution to the Serre equations with the simulated and measured results shown for ($a$) $x =3$m, ($b$) $x =4$m, ($c$) $x =4.55$m, ($d$) $x =5$m, ($e$) $x =6$m, ($f$) $x =8$m, and ($g$) $x =10.8$m using $\Delta x = 0.01115$m and $Cr = 0.2$.}
 \label{fig:Undular_Bore_Boussq_2}
\end{figure}

The recorded water surface profile at the acoustic displacement meters over time are shown in Figure \ref{fig:Undular_Bore_Boussq_2} along with the simulated water surface profile predicted by the model.

At all locations in the flume, the dispersive waves are symmetrical about the predicted water level. The simulated results show that the simulated bore speed is slightly slower than the observed bore speed. This is the theoretical observation, where the group and phase speed of waves for the Serre equations are slower than for the shallow water wave equations. These results show that the third-order scheme has accurately predicted the arrival of the bore.  In addition, it has accurately predicted the amplitude of the dispersive waves which have a slightly longer wavelength than the actual dispersive waves.

\subsubsection{Rectangular Initial Wave}

In the laboratory experiment conducted by  Hammack and Segur\cite{Hammack-Segur-1978-337},  a wave maker consists of a rectangular piston $61$cm in length at the end of a wave tank spans the full width of the tank. The tank is $31.6$m in length, $61$cm deep and $39.4$cm wide, horizontal with vertical sides and is constructed from glass. The piston moved monotonically from its initial position, which is flush with the tank bed to its final elevation.   It can be displaced vertically up or down. The upstream wall of the wave tank adjacent to the wave maker is a plane of symmetry. The length of the piston, $b = 61$cm represents the half-length of a hypothetical piston occupying the region $-b < x < b$. The symmetrical problem is simulated using the numerical schemes. A rectangular wave propagates following a sudden downward $1$cm ($h_0 = 9$cm) movement of the piston. The quiescent water depth, $h_1$ is fixed at $10$cm. The water elevation is recorded at the fixed locations; $x/h_1 = 0$, $x/h_1 = 50$, $x/h_1 = 100$, $x/h_1 = 150$, and $x/h_1 = 200$, where $x/h_1 = 0$ is the downstream edge of the piston.

The upstream and downstream boundary conditions remain constant at; $h_1 = 10$cm and $u_1 = 0$m/s. In the simulations $\Delta x = 0.0005$m, $Cr = 0.2$, $\Delta t = Cr \Delta x/\sqrt{0.1g}$ and the solution is terminated at $t = 50$s.

\begin{figure}[htb]
\centering
\begin{tabular}{cc}
\begin{psfrags}%
\psfragscanon%
%
\psfrag{s02}[b][b]{\setlength{\tabcolsep}{0pt}}%
\psfrag{s03}[t][t]{\setlength{\tabcolsep}{0pt}\begin{tabular}{c}{\Large $t\sqrt{\dfrac{g}{h_1}} - \dfrac{x}{h_1}$}\end{tabular}}%
\psfrag{s04}[b][b]{\setlength{\tabcolsep}{0pt}\begin{tabular}{c}{\Large $\dfrac{h - h_1}{h_1}$}\end{tabular}}%
%
\psfrag{x01}[t][t]{0}%
\psfrag{x02}[t][t]{50}%
\psfrag{x03}[t][t]{100}%
\psfrag{x04}[t][t]{150}%
\psfrag{x05}[t][t]{200}%
\psfrag{x06}[t][t]{250}%
%
\psfrag{v01}[r][r]{-0.3}%
\psfrag{v02}[r][r]{}%
\psfrag{v03}[r][r]{-0.2}%
\psfrag{v04}[r][r]{}%
\psfrag{v05}[r][r]{-0.1}%
\psfrag{v06}[r][r]{}%
\psfrag{v07}[r][r]{0.0}%
\psfrag{v08}[r][r]{}%
\psfrag{v09}[r][r]{0.1}%
\psfrag{v10}[r][r]{}%
\psfrag{v11}[r][r]{0.2}%
%
\resizebox{6cm}{!}{\includegraphics{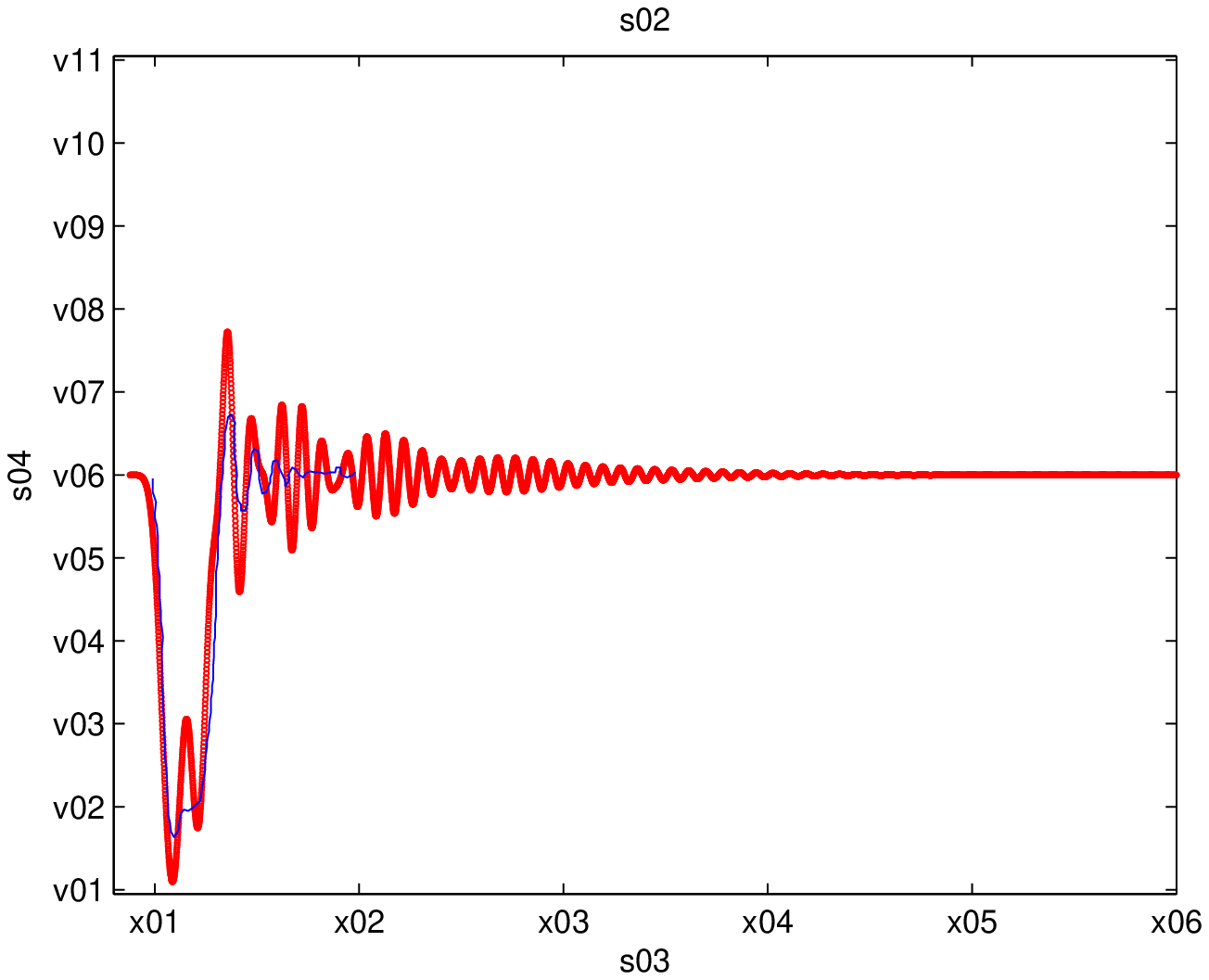}}%
\end{psfrags}%
& %
\begin{psfrags}%
\psfragscanon%
%
\psfrag{s02}[b][b]{\setlength{\tabcolsep}{0pt}}%
\psfrag{s03}[t][t]{\setlength{\tabcolsep}{0pt}\begin{tabular}{c}{\Large $t\sqrt{\dfrac{g}{h_1}} - \dfrac{x}{h_1}$}\end{tabular}}%
\psfrag{s04}[b][b]{\setlength{\tabcolsep}{0pt}\begin{tabular}{c}{\Large $\dfrac{h - h_1}{h_1}$}\end{tabular}}%
%
\psfrag{x01}[t][t]{0}%
\psfrag{x02}[t][t]{50}%
\psfrag{x03}[t][t]{100}%
\psfrag{x04}[t][t]{150}%
\psfrag{x05}[t][t]{200}%
\psfrag{x06}[t][t]{250}%
%
\psfrag{v01}[r][r]{-0.3}%
\psfrag{v02}[r][r]{}%
\psfrag{v03}[r][r]{-0.2}%
\psfrag{v04}[r][r]{}%
\psfrag{v05}[r][r]{-0.1}%
\psfrag{v06}[r][r]{}%
\psfrag{v07}[r][r]{0.0}%
\psfrag{v08}[r][r]{}%
\psfrag{v09}[r][r]{0.1}%
\psfrag{v10}[r][r]{}%
\psfrag{v11}[r][r]{0.2}%
%
\resizebox{6cm}{!}{\includegraphics{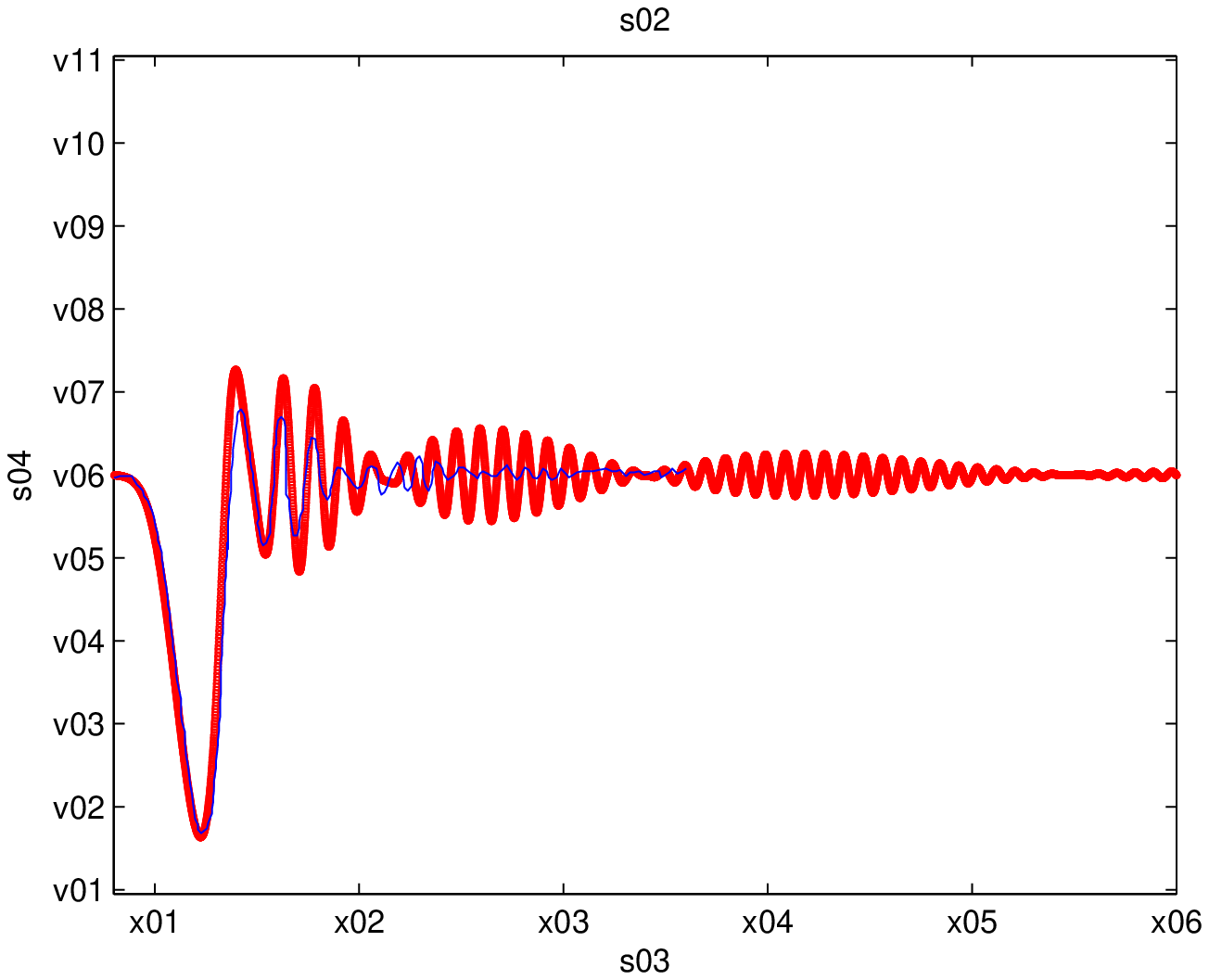}}%
\end{psfrags}%
\\
\phantom{x} & \\
($a$) & ($b$) \\
\phantom{x} & \\
\begin{psfrags}%
\psfragscanon%
%
\psfrag{s02}[b][b]{\setlength{\tabcolsep}{0pt}}%
\psfrag{s03}[t][t]{\setlength{\tabcolsep}{0pt}\begin{tabular}{c}{\Large $t\sqrt{\dfrac{g}{h_1}} - \dfrac{x}{h_1}$}\end{tabular}}%
\psfrag{s04}[b][b]{\setlength{\tabcolsep}{0pt}\begin{tabular}{c}{\Large $\dfrac{h - h_1}{h_1}$}\end{tabular}}%
%
\psfrag{x01}[t][t]{0}%
\psfrag{x02}[t][t]{50}%
\psfrag{x03}[t][t]{100}%
\psfrag{x04}[t][t]{150}%
\psfrag{x05}[t][t]{200}%
\psfrag{x06}[t][t]{250}%
%
\psfrag{v01}[r][r]{-0.3}%
\psfrag{v02}[r][r]{}%
\psfrag{v03}[r][r]{-0.2}%
\psfrag{v04}[r][r]{}%
\psfrag{v05}[r][r]{-0.1}%
\psfrag{v06}[r][r]{}%
\psfrag{v07}[r][r]{0.0}%
\psfrag{v08}[r][r]{}%
\psfrag{v09}[r][r]{0.1}%
\psfrag{v10}[r][r]{}%
\psfrag{v11}[r][r]{0.2}%
%
\resizebox{6cm}{!}{\includegraphics{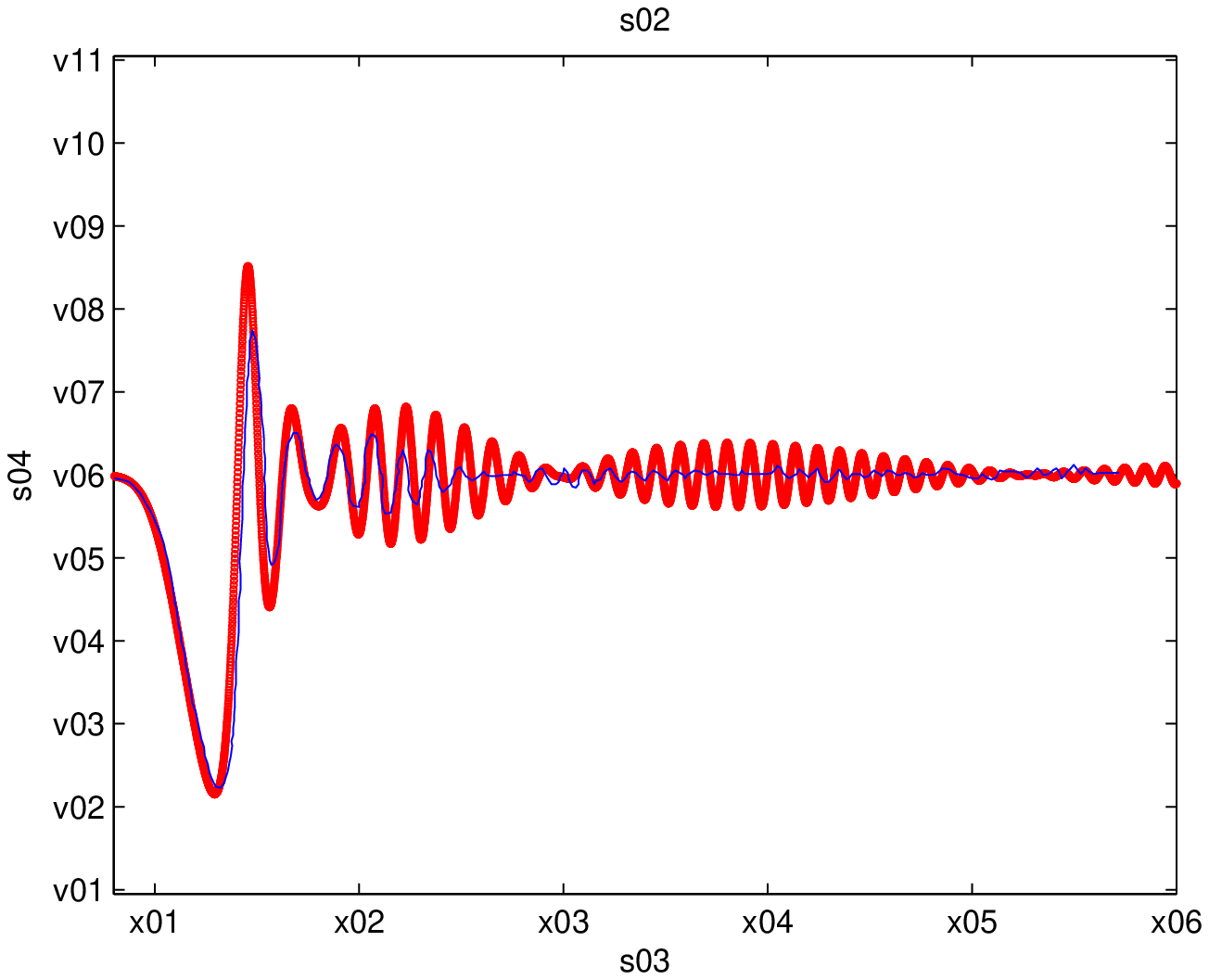}}%
\end{psfrags}%
&
\begin{psfrags}%
\psfragscanon%
%
\psfrag{s02}[b][b]{\setlength{\tabcolsep}{0pt}}%
\psfrag{s03}[t][t]{\setlength{\tabcolsep}{0pt}\begin{tabular}{c}{\Large $t\sqrt{\dfrac{g}{h_1}} - \dfrac{x}{h_1}$}\end{tabular}}%
\psfrag{s04}[b][b]{\setlength{\tabcolsep}{0pt}\begin{tabular}{c}{\Large $\dfrac{h - h_1}{h_1}$}\end{tabular}}%
%
\psfrag{x01}[t][t]{0}%
\psfrag{x02}[t][t]{50}%
\psfrag{x03}[t][t]{100}%
\psfrag{x04}[t][t]{150}%
\psfrag{x05}[t][t]{200}%
\psfrag{x06}[t][t]{250}%
%
\psfrag{v01}[r][r]{-0.3}%
\psfrag{v02}[r][r]{}%
\psfrag{v03}[r][r]{-0.2}%
\psfrag{v04}[r][r]{}%
\psfrag{v05}[r][r]{-0.1}%
\psfrag{v06}[r][r]{}%
\psfrag{v07}[r][r]{0.0}%
\psfrag{v08}[r][r]{}%
\psfrag{v09}[r][r]{0.1}%
\psfrag{v10}[r][r]{}%
\psfrag{v11}[r][r]{0.2}%
%
\resizebox{6cm}{!}{\includegraphics{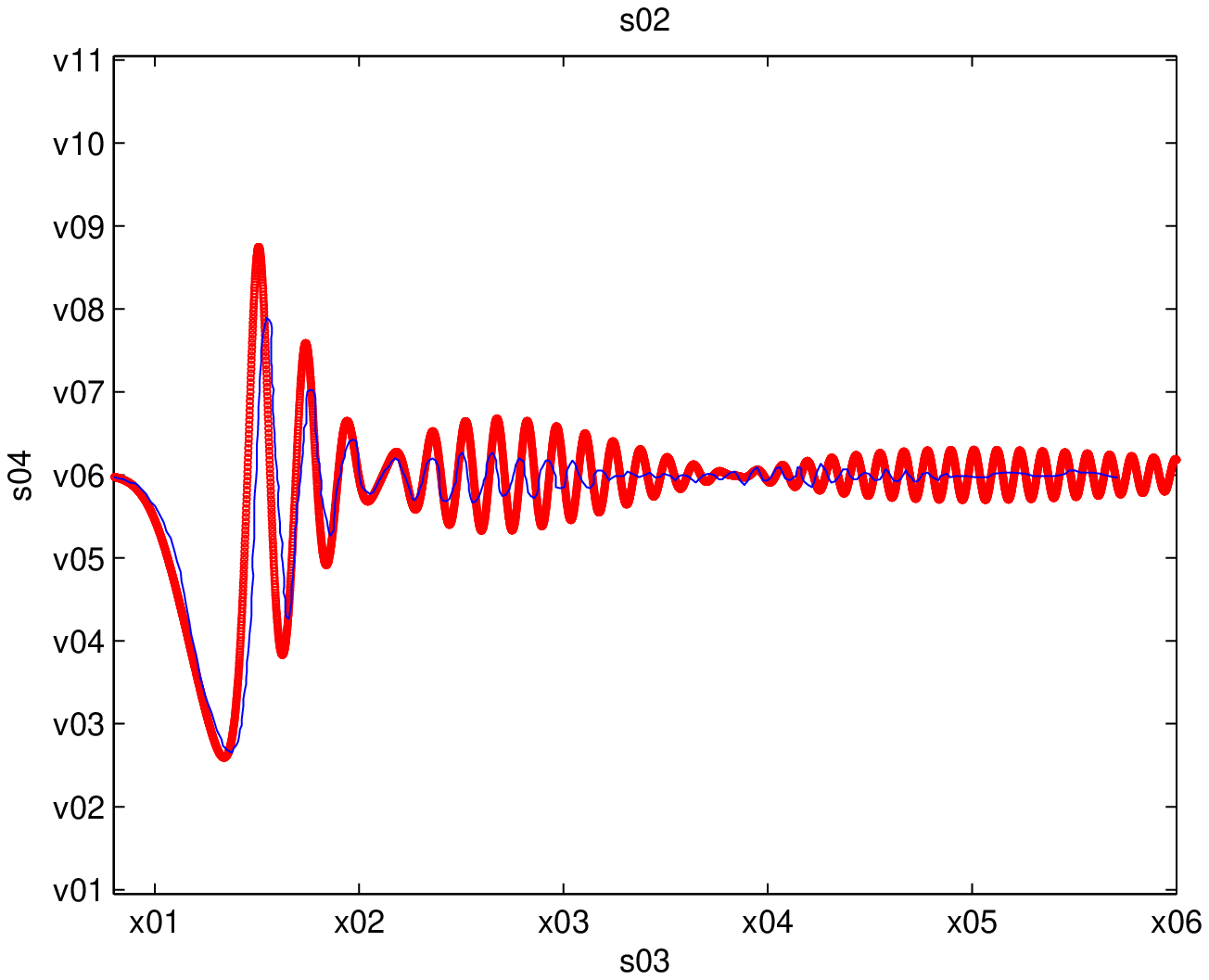}}%
\end{psfrags}%
\\
\phantom{x} & \\
($c$) & ($d$) \\
\phantom{x} & \\
\multicolumn{2}{c}{
\begin{psfrags}%
\psfragscanon%
%
\psfrag{s02}[b][b]{\setlength{\tabcolsep}{0pt}}%
\psfrag{s03}[t][t]{\setlength{\tabcolsep}{0pt}\begin{tabular}{c}{\Large $t\sqrt{\dfrac{g}{h_1}} - \dfrac{x}{h_1}$}\end{tabular}}%
\psfrag{s04}[b][b]{\setlength{\tabcolsep}{0pt}\begin{tabular}{c}{\Large $\dfrac{h - h_1}{h_1}$}\end{tabular}}%
%
\psfrag{x01}[t][t]{0}%
\psfrag{x02}[t][t]{50}%
\psfrag{x03}[t][t]{100}%
\psfrag{x04}[t][t]{150}%
\psfrag{x05}[t][t]{200}%
\psfrag{x06}[t][t]{250}%
%
\psfrag{v01}[r][r]{-0.3}%
\psfrag{v02}[r][r]{}%
\psfrag{v03}[r][r]{-0.2}%
\psfrag{v04}[r][r]{}%
\psfrag{v05}[r][r]{-0.1}%
\psfrag{v06}[r][r]{}%
\psfrag{v07}[r][r]{0.0}%
\psfrag{v08}[r][r]{}%
\psfrag{v09}[r][r]{0.1}%
\psfrag{v10}[r][r]{}%
\psfrag{v11}[r][r]{0.2}%
%
\resizebox{6cm}{!}{\includegraphics{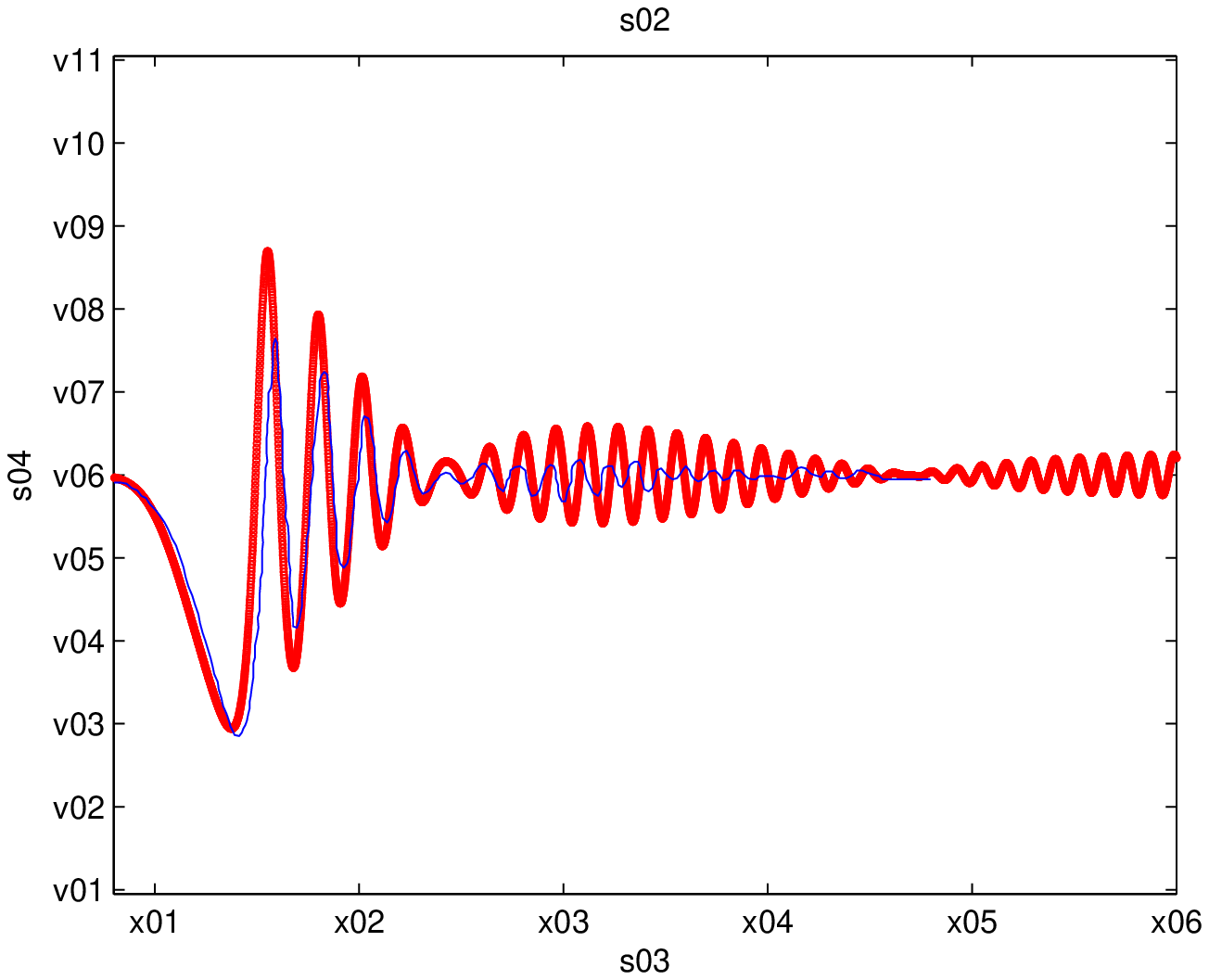}}%
\end{psfrags}%
}
\\
\phantom{x} & \\
\multicolumn{2}{c}{$(e)$} \\
\phantom{x} & \\
\end{tabular}
\caption{Measured (\textemdash) and simulated ($\circ$) water depth, $h(x,t)$ for the rectangular wave experiment in a frictionless rectangular channel,  with $h_1 = 0.1$m, $u_1 = u_0 = 0$m/s and  $h_0 = 0.09$m using second-order central and SSP Runge-Kutta scheme solution to the Serre equations with the simulated and measured results shown for the simulated and measured results shown for  ($a$) $x/h_0 = 0$, ($b$) $x /h_0 =50$, ($c$) $x/h_0 =100$, ($d$) $x/h_0 =150$,  and ($e$) $x/h_0 =200$ with $\Delta x = 0.005$m, $\Delta t = Cr \Delta x/\sqrt{gh_0}$ and $Cr = 0.2$.}
 \label{fig:Segur_Figure3_second}
\end{figure}

There are several wave train following the dominant first wave train, see Figure \ref{fig:Segur_Figure3_second}. The amplitude of the dispersive waves in these is much smaller than the amplitude of the dispersive waves that immediately follow the shock.

Results from the third-order scheme, are also shown in Figure \ref{fig:Segur_Figure3_second}. There is excellent agreement between the simulated and observed results. The rarefaction wave, shock speed and the phase of the dispersive waves are faithfully reproduced by the numerical scheme. In addition, secondary wave trains are also reproduced by the numerical scheme, see for example Figure \ref{fig:Segur_Figure3_second}($c$) at $t\sqrt{g/h_1} - x/h_1 \approx 110$. The amplitude of the dispersive waves, though are slightly overestimated by the numerical scheme. This is not surprising since the Serre equations do not contain terms that represent the internal viscosity of the fluid.

\subsection{Dam-break}

The dam-break problem is a standard test for models used to solve the shallow water wave equations, which has a known analytical solution (see, for example Zoppou and Roberts\cite{Zoppou-Roberts-2003-11}). It has been chosen to demonstrate the flexibility of the proposed model for simulating both subcritical and supercritical problems.

The dam-break problem is solved using the proposed third-order solution to the Serre equations. The simulated results have been plotted against the analytical solution to the shallow water wave equations for the dam break problem, which is used as reference data. The dam-break occurs in a frictionless rectangular channel, $1000$m in length where the initial velocity of the water $\bar{u} = 0$m/s and the water depth upstream of the dam, which is located at $x = 500$m is given by $h_1$ and downstream of the dam by $h_0$.  In all the models,  $\Delta x = 0.1$m, $Cr = 0.2$, $\Delta t = Cr \Delta x/\sqrt{gh_1}$ and the solution is terminated at $t = 30$s. Three cases are considered; $h_1 = 10$m with $h_0 = 1$m,  $h_1 = 10$m with $h_0 = 2$m and $h_1 = 1.8$m with $h_0 = 1$m. These have as their maximum Froude numbers; $Fr = u/\sqrt{gh} = 1.18$, $0.81$ and $0.29$ respectively, which were obtained from the analytical solution to the shallow water wave equations. The three problems involve supercritical flows, near critical flow and subcritical flows.

\begin{figure}[htb]
\centering
\begin{tabular}{cc}
\begin{psfrags}%
\psfragscanon%
%
\psfrag{s03}[t][t]{\setlength{\tabcolsep}{0pt}\begin{tabular}{c}{\Large $x$(m)}\end{tabular}}%
\psfrag{s04}[b][b]{\setlength{\tabcolsep}{0pt}\begin{tabular}{c}{\Large $h$(m)}\end{tabular}}%
%
\psfrag{x01}[t][t]{0}%
\psfrag{x02}[t][t]{100}%
\psfrag{x03}[t][t]{200}%
\psfrag{x04}[t][t]{300}%
\psfrag{x05}[t][t]{400}%
\psfrag{x06}[t][t]{500}%
\psfrag{x07}[t][t]{600}%
\psfrag{x08}[t][t]{700}%
\psfrag{x09}[t][t]{800}%
\psfrag{x10}[t][t]{900}%
\psfrag{x11}[t][t]{1000}%
%
\psfrag{v01}[r][r]{0}%
\psfrag{v02}[r][r]{2}%
\psfrag{v03}[r][r]{4}%
\psfrag{v04}[r][r]{6}%
\psfrag{v05}[r][r]{8}%
\psfrag{v06}[r][r]{10}%
\psfrag{v07}[r][r]{12}%
%
\resizebox{6cm}{!}{\includegraphics{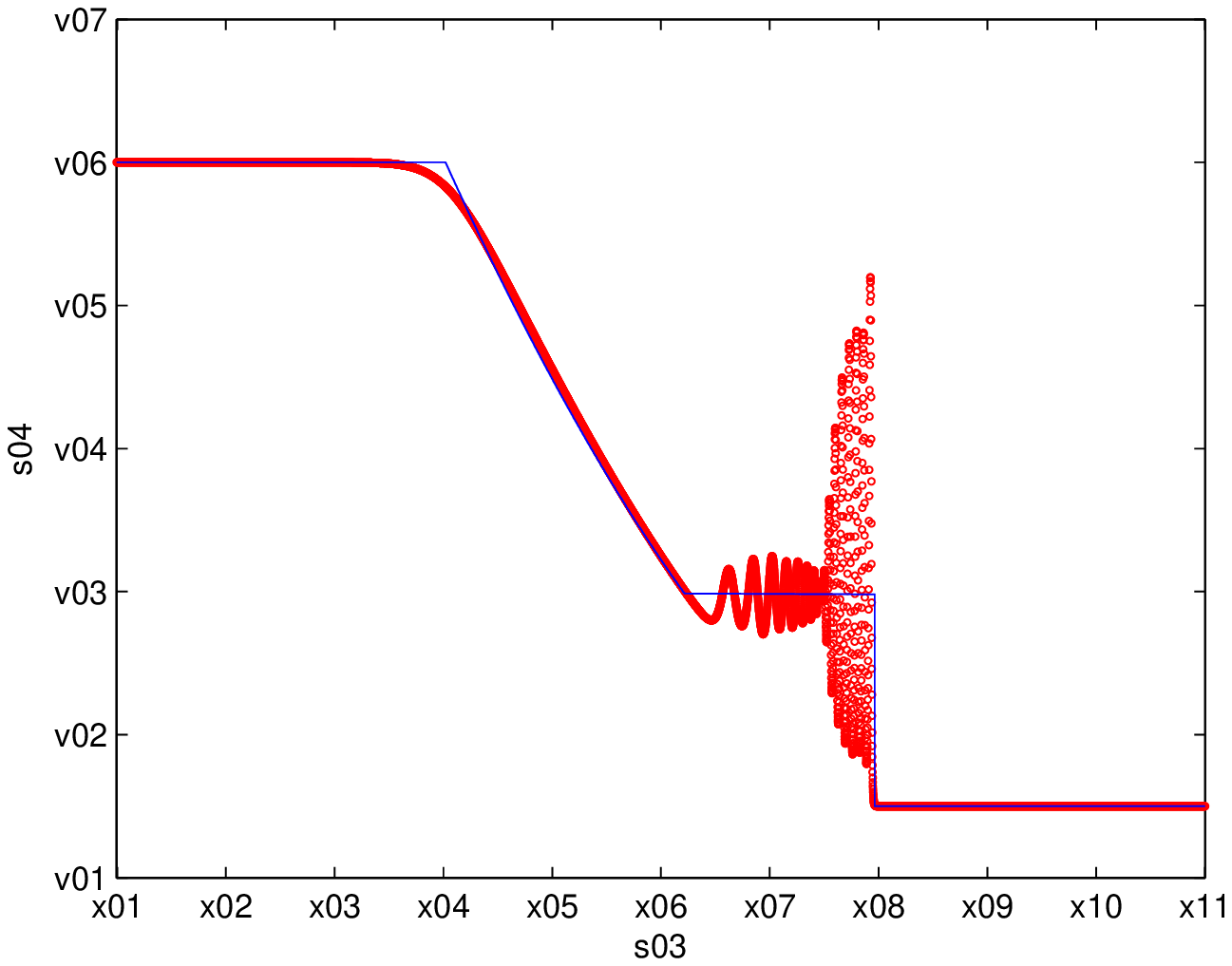}}%
\end{psfrags}%
& %
\begin{psfrags}%
\psfragscanon%
%
\psfrag{s03}[t][t]{\setlength{\tabcolsep}{0pt}\begin{tabular}{c}{\Large $x$(m)}\end{tabular}}%
\psfrag{s04}[b][b]{\setlength{\tabcolsep}{0pt}\begin{tabular}{c}{\Large $h$(m)}\end{tabular}}%
%
\psfrag{x01}[t][t]{0}%
\psfrag{x02}[t][t]{100}%
\psfrag{x03}[t][t]{200}%
\psfrag{x04}[t][t]{300}%
\psfrag{x05}[t][t]{400}%
\psfrag{x06}[t][t]{500}%
\psfrag{x07}[t][t]{600}%
\psfrag{x08}[t][t]{700}%
\psfrag{x09}[t][t]{800}%
\psfrag{x10}[t][t]{900}%
\psfrag{x11}[t][t]{1000}%
%
\psfrag{v01}[r][r]{0}%
\psfrag{v02}[r][r]{2}%
\psfrag{v03}[r][r]{4}%
\psfrag{v04}[r][r]{6}%
\psfrag{v05}[r][r]{8}%
\psfrag{v06}[r][r]{10}%
\psfrag{v07}[r][r]{12}%
%
\resizebox{6cm}{!}{\includegraphics{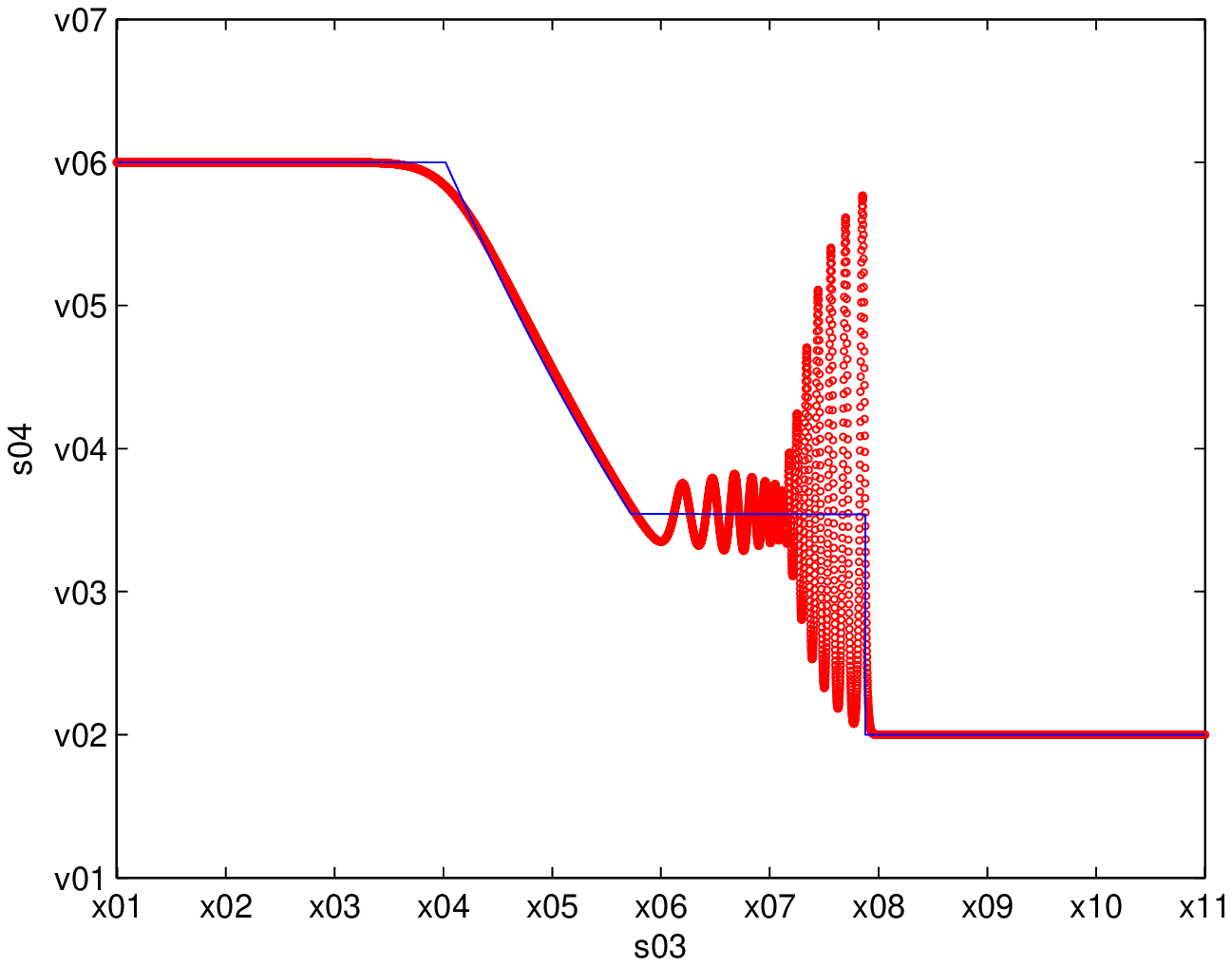}}%
\end{psfrags}%
\\
\phantom{x} & \\
($a$) & ($b$) \\
\phantom{x} & \\
\multicolumn{2}{c}{
\begin{psfrags}%
\psfragscanon%
%
\psfrag{s03}[t][t]{\setlength{\tabcolsep}{0pt}\begin{tabular}{c}{\Large $x$(m)}\end{tabular}}%
\psfrag{s04}[b][b]{\setlength{\tabcolsep}{0pt}\begin{tabular}{c}{\Large $h$(m)}\end{tabular}}%
%
\psfrag{x01}[t][t]{0}%
\psfrag{x02}[t][t]{100}%
\psfrag{x03}[t][t]{200}%
\psfrag{x04}[t][t]{300}%
\psfrag{x05}[t][t]{400}%
\psfrag{x06}[t][t]{500}%
\psfrag{x07}[t][t]{600}%
\psfrag{x08}[t][t]{700}%
\psfrag{x09}[t][t]{800}%
\psfrag{x10}[t][t]{900}%
\psfrag{x11}[t][t]{1000}%
%
\psfrag{v01}[r][r]{0.8}%
\psfrag{v02}[r][r]{1.0}%
\psfrag{v03}[r][r]{1.2}%
\psfrag{v04}[r][r]{1.4}%
\psfrag{v05}[r][r]{1.6}%
\psfrag{v06}[r][r]{1.8}%
\psfrag{v07}[r][r]{2.0}%
%
\resizebox{6cm}{!}{\includegraphics{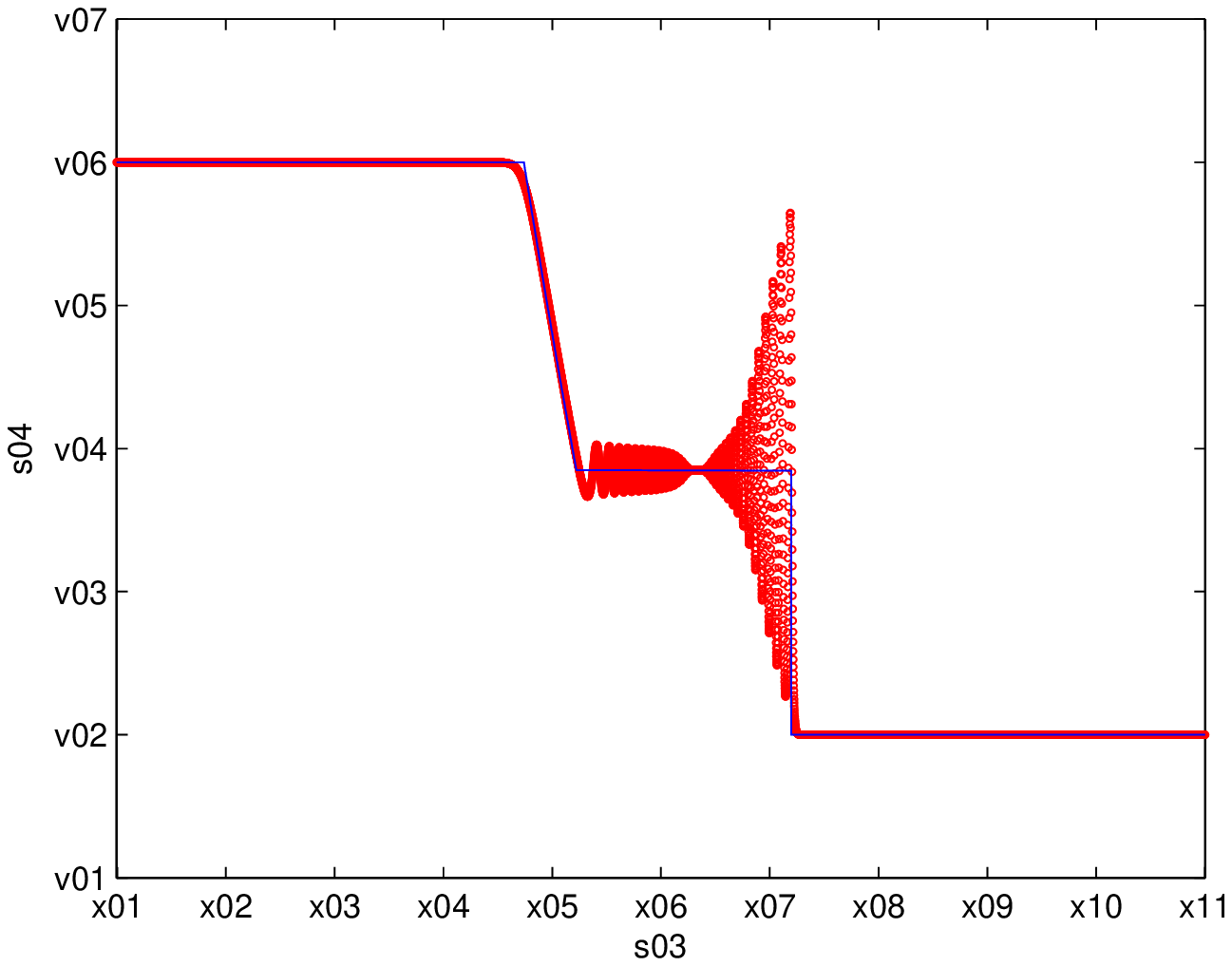}}

\end{psfrags}
}
\\
\multicolumn{2}{c}{$(c)$} \\
\phantom{x} &  \\
\end{tabular}
\caption{Analytical (\textemdash) solution to the shallow water wave equations and simulated ($\circ$) water depth, $h(x,t)$ for the dam break problem in a frictionless rectangular channel, $1000$m in length, $u_1 = u_0 = 0$m/s using the third-order Serre equations solver with ($a$) $h_1 = 10$m and $h_0 = 1$m, ($b$) $h_1 = 10$m and $h_0 = 2$m and ($c$)  $h_1 = 1.8$m and $h_0 = 1$m with $\Delta x = 0.1$m, $Cr = 0.2$, $\Delta t = Cr \Delta x/\sqrt{gh_1}$ and the solution is terminated at $t = 30$s.}
 \label{fig:Dam_break_10_1}
\end{figure}

The simulated results using are shown in Figures \ref{fig:Dam_break_10_1}. The arrival of the shock is accurately captured, as is the rarefaction fan and the shock height. The results for the third-order scheme, shown in Figure \ref{fig:Dam_break_10_1}($c$) are very similar to those obtained by El \emph{et al.}\cite{El-etal-2006}, who used a second-order Lax-Wendroff scheme and the second-order finite volume model described by Zoppou \emph{et al.}\cite{Zoppou-etal-2016} to solve the Serre equations.

An interesting feature of the results shown in Figure \ref{fig:Dam_break_10_1} is that the oscillations are bounded. They are restricted to the minimum and maximum initial water depth. The simulated water velocity is also bounded. The Serre equations also conserve energy. Since energy for the Serre equations is uniformly bounded\cite{Li-Y-2006-1255} then the solution is also bounded by the energy of the initial conditions.

\section{Conclusions}

By replacing the mix derivative term in the flux term by a combination of  temporal and spatial derivative terms, the Serre equations can be written in conservation law form, where the system of homogeneous equations contains a new conserved quantity and its corresponding flux term. The water depth and the new conserved quantity is evolved using a finite volume scheme. The water velocity, which is the remaining primitive variable is obtained by solving a second-order elliptic equation using finite elements.

Using analytical solutions, laboratory flume data and by simulating the dam-break problem, the proposed third-order hybrid finite volume/finite element scheme is shown to be simple to implement and stable for a range of problems including rapidly varying flows.  It accurately predicts the phase, arrival of the dispersive waves and their amplitude that are associated with rapidly varying flows.

The proposed scheme is currently being extended to problems with bathymetry, to two-dimensional problems and to other system of dispersive equations.

\section*{Acknowledgements}

Professor H. Chanson, Department of Civil Engineering, University of Queensland for providing the data for the undular bore and to Dr David George, Cascades Volcano Observatory, U.S. Geological Survey for providing the rectangular wave data.

\appendix

\begin{landscape}
\section{Stiffness and Load Marix in the Finite Elemen Scheme}

The coefficients in the stiffness matrix $\mathbf{Q}_e$ are given by
\begin{small}
\begin{equation*}
\left [
\begin{array}{c}
q_{11} \\
q_{12} \\
q_{13}
\end{array}
\right ] =  \dfrac{1}{3780 \Delta x}  \left [
\begin{array}{c}
912 h^+_{j-1/2} h_j^2 + 948 (h^+_{j-1/2})^2 h_j + 61 (h^-_{j+1/2})^3 + 21 h^+_{j-1/2}( h^-_{j+1/2})^2 - 195 (h^+_{j-1/2})^2 h^-_{j+1/2} - 240 h_j^2 h^-_{j+1/2} - 336 h^+_{j-1/2} h_j h^-_{j+1/2} + 832 h_j^3 + 853 (h^+_{j-1/2})^3 + 84 h_j (h^-_{j+1/2})^2 \\ \\
 - 240 h_j (h^-_{j+1/2})^2 - 284 (h^-_{j+1/2})^3 - 1104 (h^+_{j-1/2})^2 h_j + 384 h^+_{j-1/2} h_j h^-_{j+1/2} - 1076 (h^+_{j-1/2})^3 - 512 h_j^3 - 960 h^+_{j-1/2} h_j^2 + 228 (h^+_{j-1/2})^2 h^-_{j+1/2} + 12 h^+_{j-1/2} (h^-_{j+1/2})^2 + 192 h_j^2 h^-_{j+1/2} \\ \\
156 (h^+_{j-1/2})^2 h_j - 33 (h^+_{j-1/2})^2 h^-_{j+1/2} - 33 h^+_{j-1/2} (h^-_{j+1/2})^2 + 156 h_j (h^-_{j+1/2})^2 + 48 h_j^2 h^-_{j+1/2} - 48 h^+_{j-1/2} h_j h^-_{j+1/2} + 223 (h^-_{j+1/2})^3 - 320 h_j^3 + 223 (h^+_{j-1/2})^3 + 48 h^+_{j-1/2} h_j^2
\end{array}
\right ]
\end{equation*}
\\
\begin{equation*}
\left [
\begin{array}{c}
q_{21} \\
q_{22} \\
q_{23}
\end{array}
\right ] =  \dfrac{1}{3780 \Delta x}  \left [
\begin{array}{c}
-240 h_j (h^-_{j+1/2})^2 - 284 (h^-_{j+1/2})^3 - 1104 (h^+_{j-1/2})^2 h_j + 384 h^+_{j-1/2} h_j h^-_{j+1/2} - 1076 (h^+_{j-1/2})^3 - 512 h_j^3 - 960 h^+_{j-1/2} h_j^2 + 228 (h^+_{j-1/2})^2 h^-_{j+1/2} + 12 h^+_{j-1/2} (h^-_{j+1/2})^2 + 192 h_j^2 h^-_{j+1/2} \\ \\
1344 (h^+_{j-1/2})^2 h_j - 240 (h^+_{j-1/2})^2 h^-_{j+1/2} + 768 h^+_{j-1/2} h_j^2 - 240 h^+_{j-1/2} (h^-_{j+1/2})^2 + 768 h_j^2 h^-_{j+1/2} + 1344 h_j (h^-_{j+1/2})^2 + 1360 (h^+_{j-1/2})^3 + 1024 h_j^3 + 1360 (h^-_{j+1/2})^3 - 768 h^+_{j-1/2} h_j h^-_{j+1/2} \\ \\
12 (h^+_{j-1/2})^2 h^-_{j+1/2} - 240 (h^+_{j-1/2})^2 h_j - 1104 h_j (h^-_{j+1/2})^2 + 192 h^+_{j-1/2} h_j^2 - 960 h_j^2 h^-_{j+1/2} - 512 h_j^3 - 1076 (h^-_{j+1/2})^3 + 228 h^+_{j-1/2} (h^-_{j+1/2})^2 - 284 (h^+_{j-1/2})^3 + 384 h^+_{j-1/2} h_j h^-_{j+1/2}
\end{array}
\right ]
\end{equation*}
\\
\begin{equation*}
\left [
\begin{array}{c}
q_{31} \\
q_{32} \\
q_{33}
\end{array}
\right ] =  \dfrac{1}{3780 \Delta x}  \left [
\begin{array}{c}
-48 h^+_{j-1/2} h_j h^-_{j+1/2} - 320 h_j^3 - 33 h^+_{j-1/2} (h^-_{j+1/2})^2 + 48 h_j^2 h^-_{j+1/2} + 156 h_j (h^-_{j+1/2})^2 + 48 h^+_{j-1/2} h_j^2 + 156 (h^+_{j-1/2})^2 h_j - 33 (h^+_{j-1/2})^2 h^-_{j+1/2} + 223 (h^+_{j-1/2})^3 + 223 (h^-_{j+1/2})^3 \\ \\
12 (h^+_{j-1/2})^2 h^-_{j+1/2} - 240 (h^+_{j-1/2})^2 h_j - 1104 h_j (h^-_{j+1/2})^2 + 192 h^+_{j-1/2} h_j^2 - 960 h_j^2 h^-_{j+1/2} - 512 h_j^3 - 1076 (h^-_{j+1/2})^3 + 228 h^+_{j-1/2} (h^-_{j+1/2})^2 - 284 (h^+_{j-1/2})^3 + 384 h^+_{j-1/2} h_j h^-_{j+1/2} \\ \\
-336 h^+_{j-1/2} h_j h^-_{j+1/2} - 195 h^+_{j-1/2} (h^-_{j+1/2})^2 + 912 h_j^2 h^-_{j+1/2} + 948 h_j (h^-_{j+1/2})^2 + 832 h_j^3 + 853 (h^-_{j+1/2})^3 + 61 (h^+_{j-1/2})^3 + 84 (h^+_{j-1/2})^2 h_j + 21 (h^+_{j-1/2})^2 h^-_{j+1/2} - 240 h^+_{j-1/2} h_j^2
\end{array}
\right ].
\end{equation*}
\end{small}
\end{landscape}

\begin{landscape}
The stiffness matrix, $\mathbf{A}$, the unknowns, $\bar{\mathbf{u}}$ and the stiffness matrix, $\mathbf{R}$ for the finite element solution to the second-order elliptic equation, \eqref{eq:G_FE} are given by
\begin{small}
\begin{equation*}
\mathbf{A} = \left [
\begin{array}{ccccccccc}
1 & 0 &  &  &  &  &  &  &  \\
q_{21}^{(1)} +  p_{21}^{(1)} & q_{22}^{(1)} +  p_{22}^{(1)} & q_{23}^{(1)} +  p_{23}^{(1)} &  &  &  &  &  &  \\
q_{31}^{(1)} + p_{31}^{(1)} & q_{32}^{(1)} +  p_{32}^{(1)} & q_{33}^{(1)} +  p_{33}^{(1)} + q_{11}^{(2)} +  p_{11}^{(2)}  & q_{12}^{(2)} +  p_{12}^{(2)}  & q_{13}^{(2)} +  p_{13}^{(2)}  &  &  &  \\
 &   & q_{21}^{(2)} +  p_{21}^{(2)}   & q_{22}^{(2)} +  p_{22}^{(2)}   & q_{23}^{(2)} +  p_{23}^{(2)}   &  &  &  &  \\
 &  & q_{31}^{(2)} +  p_{31}^{(2)}   & q_{32}^{(2)} +  p_{32}^{(2)}   & q_{33}^{(2)} +  p_{33}^{(2)} + q_{11}^{(3)} +  p_{11}^{(3)}  &  q_{12}^{(3)} +  p_{12}^{(3)}  & q_{13}^{(3)} +  p_{13}^{(3)} &  &  \\
 &  &  &  &  q_{21}^{(3)} +  p_{21}^{(3)}  & q_{22}^{(3)} +  p_{22}^{(3)}   & q_{23}^{(3)} +  p_{23}^{(3)}  & & \\
 &  &  &  & q_{31}^{(3)} +  p_{31}^{(3)} & q_{32}^{(3)} +  p_{32}^{(3)} & q_{33}^{(3)} +  p_{33}^{(3)} + q_{11}^{(4)} +  p_{11}^{(4)}  & q_{12}^{(4)} +  p_{12}^{(4)}  & q_{13}^{(4)} +  p_{13}^{(4)} \\
 &  &  &  &  &  & q_{21}^{(4)} +  p_{21}^{(4)} & q_{22}^{(4)} +  p_{22}^{(4)}  & q_{23}^{(4)} +  p_{23}^{(4)} \\
 &  &  &  &  &  &  & 0 & 1 \\
\end{array}
\right ],
\end{equation*}
\begin{equation*}
\bar{\mathbf{u}} = \left [
\begin{array}{c}
u_{j-1/2} \\
u_j \\
u_{j+1/2} \\
u_{j+1} \\
u_{j+3/2} \\
u_{j+2} \\
u_{j+5/2} \\
u_{j+3} \\
u_{j+7/2}
\end{array}
\right ]  \quad \text{and} \quad \mathbf{R} = \left [
\begin{array}{c}
\alpha \\
r_2^{(1)} \\
r_3^{(1)} + r_1^{(2)} \\
r_2^{(2)} \\
r_3^{(2)} + r_1^{(3)} \\
r_2^{(3)} \\
r_3^{(3)} + r_1^{(4)} \\
r_2^{(4)} \\
\beta
\end{array}
\right ]
\end{equation*}
\end{small}
where the superscript denotes the element number and the coefficients $\alpha$ and $\beta$ are provided by the essential Dirichlet boundary conditions.
\end{landscape}


\end{document}